\input amstex
\input amsppt.sty
\NoBlackBoxes
\nologo
\topmatter
\title A note on Toroidalization:\linebreak
the problem of resolution of singularities \linebreak
of morphisms in the logarithmic category  
\endtitle
\author Kenji Matsuki
\endauthor
\leftheadtext{Kenji Matsuki}
\rightheadtext{Toroidalization}
\endtopmatter

\document

$$\bold{Contents}$$

\S 0. Introduction

\S 1. Formulation of toroidalization as resolution of singularities of morphisms in 

\hskip.24in the
logarithmic category

\S 2. Algorithm for toroidalization

\S 3. Local analysis in the case where $\dim X = \dim Y = 2$

\S 4. Termination of the alogorithm in the case where $\dim X = \dim Y = 2$

\vskip.2in

\S 0. Introduction

\vskip.1in

0.1. $\bold{Prologue}$.  The purpose of this note is to discuss the following problem of
``toroidalization": Given a (proper) morphism
$f:X
\rightarrow Y$ between nonsingular varieties, can we find sequences of blowups with
smooth centers $\pi_X:X' \rightarrow X$ and $\pi_Y:Y' \rightarrow Y$ such that the
induced map $f':X' \rightarrow Y'$ is a toroidal morphism (See \S 1 for the precise
definition of a morphism to be toroidal.) ?
$$\CD
X \hskip.1in@. \overset{\pi_X}\to{\longleftarrow} @. \hskip.1in X' \\
@V{f}VV @. @VV{f'}V \\
Y \hskip.1in @. \overset{\pi_Y}\to{\longleftarrow} @. \hskip.1in Y' \\
\endCD$$

The problem of toroidalization evolved along with the problem of semi-stable reduction and
with the problem of factorization of birational maps, among others.  The motivation behind
is that, via toroidalization, we should be able to reduce the original problem to a
combinatorial one, utilizing the correspondence between the geometry of toric varieties
(toroidal embeddings) and the combinatorics of convex bodies.  As such, it may seem as
though the problem of toroidalization is a mere working hypothesis.  The theme of this note
is, however, to recognize its intrinsic importance, by presenting and reformulating it as the
problem of resolution of singularities of morphisms in the logarithmic category.  (See \S 1
for the precise meaning.)  This view point was first communicated to the author by Prof. Dan
Abramovich, though apparently it had been known to the experts for a while (cf.
Kato[20]Abramovich-Karu[4]).  

\vskip.1in

0.2. $\bold{A\ naive\ question}$.  A naive question can be posed in the
following way: Given a (proper and surjective) morphism $f:X \rightarrow
Y$ between nonsingular varieties, by choosing suitable ``modifications" $\pi_X:X'
\rightarrow X$ and $\pi_Y:Y' \rightarrow Y$, how ``nice" can one make the induced morphism
$f':X'
\rightarrow Y'$ ?  (We require that we choose modifications $\pi_X$ and $\pi_Y$ in such
a way that $f'$ is a morphism, which a priori is just a rational map.)  Of course, in
order to make this question meaningful, one has to make explicit the mathematical
specifications of the words ``modifications" and ``nice".

When $\dim Y = 1$ and we specify the morphism $f'$ to be ``nice" if every fiber is
reduced with simple normal crossings, while the modification $\pi_Y$ is restricted to
a finite morphism (between nonsingular varieties) and the modification $\pi_X$ is
restricted to a sequence of blowups with smooth centers (possibly via
normalization after base change), the question becomes the problem of semi-stable
reduction (over a base of dimension one).

When $\dim Y = 1$ and we specify the morphism $f'$ to be ``nice" if every fiber has
only simple normal crossings, while the modification $\pi_Y$ is restricted to an
isomorphism and the modification $\pi_X$ is restricted to a sequence of blowups with
smooth centers, the question becomes the problem of (embedded) resolution of
hypersurface singularities.

When $f$ is birational and we specify the morphism $f'$ to be ``nice" if it is an
isomorphism, while the modifications $\pi_X$ and $\pi_Y$ are restricted to sequences
of blowups with smooth centers, the question is the strong factorization conjecture of
birational maps.

It should be warned, however, that the way we pose our question is very naive: for
certain choices of the specifications of the word ``modifications" or ``nice", we often
obtain a question to which the answer is \it no \rm in general.  This can be observed,
for example, when we specify the morphism $f'$ to be ``nice" if it is flat, while the
modifications
$\pi_X$ and $\pi_Y$ are restricted to sequences of blowups with smooth centers.  This
flattening process by blowups with smooth centers does not exist in general, even for a
generically finite morphism between nonsingular surfaces $f:X \rightarrow Y$, as is
demonstrated by a famous example by Abhyankar[1] for the failure of simultaneous
resolution of singularities.  (The reader is also invited to see an easy toric example
by T. Katsura in Oda[27].)  It is, meanwhile, an interesting question whether the above
flattening process by blowups with smooth centers exists or not, under the extra
condition that the original $f$ is a morphism with connected fibers.  (The author
does not know whether the answer is affirmative or negative.)

\vskip.1in 
 
The problem
of toroidalization appears as one form of the questions in the above setting,
where we specify $f'$ to be ``nice" if
$f'$ is toroidal, while the modifications
$\pi_X$ and $\pi_Y$ are restricted to sequences of blowups with smooth centers.

\vskip.1in

It should be mentioned and emphasized that Abramovich-Karu[4] succeeds in making the
morphism $f'$ toroidal, starting from a given proper morphism $f$, but by \it not
\rm restricting the modifications $\pi_X$ and $\pi_Y$ to be only sequences of blowups
with smooth centers.

\vskip.1in

0.3. $\bold{Algorithm}$. The reformulation of toroidalization as log resolution,
though simple and straightforward as it is, naturally leads one to the following
algorithmic approach (See \S 2 for details.): 

Step 1: Characterize the locus where the morphism is not toroidal as the logarithmic
ramification locus, which can then be computed using the logarithmic differential forms.

Step 2: Choose the center of blowup downstairs on $Y$, which should be contained in the
image of the logarithmic ramification locus under $f$.

Step 3: Apply an algorithm for canonical principalization to the pull-back upstairs on
$X$ of the defining ideal of the center chosen in Step 2 downstairs on $Y$.

If the induced morphism is toroidal, then we stop.  If not, we go back to Step 1
with the original morphism replaced by the induced one.

Termination: Show that the algorithm terminates after finitely many steps.

In fact, in the case where $\dim X = \dim Y = 2$, there is essentially no choice involved in
the above algorithm and hence is uniquely determined.  Thus the only remaining problem is to
show termination.  We accomplish this goal by providing a local analysis of the algorithm
and then showing that the logarithmic ramification locus \it decreases \rm (under certain
order, with respect to which the descending chain condition is satisfied).

\vskip.1in

A discussion of the higher dimensional case will be published elsewhere.

\vskip.1in

0.4. $\bold{A\ possible\ approach\ to\ the\ factorization\ problem}$.  The problem of
\linebreak
toroidalization was conceived as a possible approach to the (strong)
factorization conjecture of birational maps (between nonsingular varieties into a
sequence of blowups and blowdowns with smooth centers) in two steps:

Step I. Given a proper birational morphism $f:X \rightarrow Y$ (if you start with a
birational map $f:X \dashrightarrow Y$, then use the elimination of indeterminacy by
Hironaka to reach the stage where $f$ can be assumed to be a morphism), apply
toroidalization to make $f':X' \rightarrow Y'$ toroidal.

Step II. Factor the toroidal birational morphism $f':X' \rightarrow Y'$ into a sequence
consisting of blowups with smooth centers immediately followed by blowdowns with
smooth centers, by applying the toroidal version of a combinatorial
algorithm for strong factorization of toric birational maps.

(Today the strong factorization conjecture remains open, while the weak factorization
conjecture, which is strong enough for most of the applications, is a theorem
achieving Steps I and II in a slightly weaker sense:

Step I$'$: Given a proper birational morphism $f:X \rightarrow Y$, apply the theory of
birational cobordism by W{\l}odarczyk[29] and torification by
Abramovich-Karu-Matsuki-W{\l}odarczyk[6] to decompose $f$ as a composite of toroidal
birational maps $f_i:X_i \dashrightarrow Y_i \text{\ for\ }i = 0, ... , l$ (where
$X = X_0, Y_i = X_{i+1} \text{\ for\ }i = 0, ... , l-1$, and $Y_l = Y$).

Step II$'$: Apply the toroidal version of the algorithm for weak factorization by Morelli
(cf. Abramovich-Matsuki-Rashid[5]Matsuki[23]W{\l}odarczyk[28]) to each toroidal birational
map
$f_i:X_i
\dashrightarrow Y_i$.

We fall short of realizing strong factorization, since in Step I$'$ the toroidal
structures for $f_i$'s are not compatible with each other in general, and since in Step
II$'$ we know only of the algorithm for weak (but not strong) factorization even for
toric birational maps (in dimension three or more) for the moment.)

\vskip.1in

0.5. $\bold{Work\ of\ S.D. Cutkosky}$.  Even at the time when, inspired by the work of
A.J. de Jong, we could only daydream of the above approach to the problem of
factorization, S.D. Cutkosky[10] made a remarkable announcement of a proof for \it local
\rm factorization of birational maps in dimension 3.  (The \it local \rm version of the
problem of factorization replaces a birational map between nonsingular varieties $X$
and $Y$ with birational local rings dominated by a fixed valuation, and replaces
blowups along smooth centers with monoidal transforms followed by localization
subordinate to the valuation.)  After Cutkosky gave a series of lectures on his
breakthrough at Purdue University, we communicated our approach to \it global \rm
factorization via Steps I and II as above, suggesting that the
\it local
\rm version of toroidalization may be called ``monomialization" and that the first
non-trivial and plausible test case of toroidalization should be checked when $\dim X =
3$ and $\dim Y = 2$.

In a sequence of papers [9][10][11], Cutkosky brilliantly establishes the monomialization in
arbitrary dimension, thus completing Step I in the \it local \rm case.

Recently Karu[19] succeeded in providing an algorithm for \it local \rm strong factorization
of toric (toroidal) birational maps in arbitrary dimension, extending the earlier work of
Christensen[9] in dimension 3.

Thus the \it local \rm strong factorization of birational maps, known also as the Abhyankar
conjecture, is now a theorem according to the general scheme of establishing Steps I and II.

Moreover, Cutkosky[12] also establishes the \it global \rm version, i.e., the toroidalization
when $\dim X = 3$ and $\dim Y = 2$.

\vskip.1in

0.6. $\bold{Contents\ of\ our\ paper}$.  The simple purpose of this paper is to give
a proof of toroidalization in dimension 2, guided by the
principle \linebreak
$\bold{toroidalization = log\ resolution}$, which has been the main point of our theme
from the conception of the problem (cf.
Abramovich-Karu-Matsuki-W{\l}odarczyk[6]\linebreak
Matsuki[24]).  It is the belief (prejudice) of the author that only thorugh this
reformulation the ultimate solution would be achieved by finding a much-sought-for
inductional structure needed to prove the general case of toroidalization.

\vskip.1in

The contents of the paper are as follows.

\vskip.1in

In \S 1, we present our formulation of the problem of toroidalization as the problem
of resolution of singularities of morphisms in the logarithmic category.  The key is
the study of the behavior of the logarithmic differential forms.  In \S 2, we present an
essentially unique algorithm for toroidalization in the case where $\dim X = \dim Y = 2$.  In
\S 3, we present a proof for toroidalization in dimension 2.  Several proofs are known
for toroidalization in dimension 2 by now.  However, all the existing ones use the strong
factorization theorem of a (proper) birational \it morphism \rm (via Castelnuovo's
contractibility criterion of a $(-1)$-curve) in the process of verification.  This has a
fatal defect: not only one cannot hope to extend the process to higher dimension where the
factorization of a proper birational \it morphism \rm into blowups with smooth centers fails
to hold, but also it would bring a logical loop in our approach to solving the
problem of factorization through toroidalization, even in dimension 2.  They also use the
properness assumption of
$f$ in an essential way, and as a result, their analysis is local on $Y$ but global on $X$
in nature.  In contrast, we provide a detailed but rather elementary analysis of the
algorithm, which is completely local both on
$Y$ and on $X$.  As a consequence, we show
toroidalization in dimension 2 without the properness assumption of $f$.  We do not use the
strong factorization theorem in the process.  Therefore, establishing Steps I and II of the
general scheme, we give a new proof of the strong factorization of birational maps
in diemsnion 2.

The way we show that the algorithm for toroidalization teminates is also guided by the
principle, rather than setting up invariants in an adhoc way.  We show that
the logarithmic ramification locus \it decreases \rm (under certain order, with respect to
which the descending chain condition is satisfied) in the process.  

\vskip.1in

0.7. $\bold{Assumption\ on\ the\ base\ field}$.  In this paper, we assume that the base
field
$k$ is an algebraically closed field of characteristic zero.  The assumption of $k$
being algebraically closed is purely for the simplicity of the presentation: the
``canonical" nature of our algorithm (and of the algorithm for principalization of ideals)
implies that the process prescribed over
$\overline{k}$ is equivariant under the action of the Galois group
$\roman{Gal}(\overline{k}/k)$ even when $k$ is not algebraically closed, and hence that
it is actually defined over $k$.  Therefore, this assumption can be removed without much
further ado.  The assumption of $k$ being of characteristic zero, however, is essential,
because the problem of toroidalization as we formulate is too restrictive in positive
characteristic (Some easy counter examples to this restrictive formulation can be found in
positive characteristic.  See, e.g. Cutkosky-Piltant[13].) and hence needs modification.

\vskip.1in

0.7. $\bold{Acknowledgements}$. \hskip.1in The author would like to thank Dan Abramovich for
guidance and help in the subject.  Thanks are also due to Kalle Karu and Jaros{\l}aw
W{\l}odarczyk for many invaluable comments.  The author would like to thank the organizers of
the Barrett conference held at the University of Tennessee at Knoxville in April 2002,
especially Yasuyuki Kachi, Shashikant Mulay and Pavlos Tzermias, for providing a
mathematically stimulating atmosphere through the excellent lecture series by Yujiro Kawamata
and Yasutaka Ihara.

\newpage

\S 1. Formulation of toroidalization as resolution of singularities of morphisms in the
logarithmic category

\vskip.1in

1.1. $\bold{Outline\ of\ this\ section}$. \hskip.1in In this section, we formulate the
problem of toroidalization and the problem of resolution of singularities of morphisms
in the logarithmic category, and show that a solution to the latter would imply a
solution to the former.  (It could be said that, modulo some minor technical points, the
two problems are essentially equivalent.)

For this purpose, we define when a morphism $f:X \rightarrow Y$ between nonsingular
varieties is toroidal, or in the presence of the prescribed structures as
nonsingular toroidal embeddings, when a morphism $f:(U_X,X) \rightarrow (U_Y,Y)$
between nonsingular toroidal embeddings is toroidal (with respect to the given toroidal
structures
$(U_X,X)$ and
$(U_Y,Y)$).  We also define, in the presence of the logarithmic structures as
nonsingular toroidal embeddings, when a morphism
$f:(U_X,X)
\rightarrow (U_Y,Y)$ in the logarithmic category is (log) smooth.  (In order to make a
clear distinction between the notion of a morphism being ``smooth" in the usual category
and the notion of a morphism being ``smooth" in the logarithmic category, we always use
the word ``log smooth" for the latter.)  One of the main points of the section is,
though it is quite straightforward, to show that the condition of $f$ being toroidal is
equivalent to $f$ being log smooth.

\vskip.1in

1.2. $\bold{Basic\ definitions}.$ 

\vskip.1in

1.2.1. \it Logarithmic category. \rm \hskip.1in An
object in the logarithmic category (of nonsingular toroidal embeddings) is a
nonsingular toroidal embedding
$(U_X,X)$, i.e., a pair consisting of a nonsingular variety $X$ and a dense open subset
$U_X$ such that
$D_X = X - U_X$ is a divisor with only simple normal crossings.  

(A divisor $D_X$ with only simple normal crossings consists of
\it smooth \rm irreducible components and satisfies the condition that at any closed
point
$p
\in D_X$ there exists an open analytic neighborhood $p \in U_p \subset X$ with a system
of regular parameters $(x_1, ... , x_d)$ of $\widehat{{\Cal O}_{X,p}}$ such that
$$D_X \cap U_p = \{\prod_{m \in M}x_m = 0\}$$
for some subset $M \subset \{1, ... , d = \dim X\}$.)

A morphism $f:(U_X,X) \rightarrow (U_Y,Y)$ in the logarithmic category is a morphism
between nonsingular toroidal embeddings such that:

(a) $f:X \rightarrow Y$ is a dominant morphism between nonsingular varieties,

(b) $D_X = f^{-1}(D_Y)$,

(c) $f|_{U_X}:U_X \rightarrow U_Y$ is a smooth morphism.

\vskip.1in

1.2.2. \it Definition of a morphism being toroidal. \rm \hskip.1in Let $f:(U_X,X)
\rightarrow (U_Y,Y)$ be a morphism in the logarithmic category of nonsingular toroidal
embeddings.  We say $f$ is toroidal if for any $p \in X$ with $q = f(p) \in Y$ there
exist a toric (equivariant) morphism $\varphi:(T_V,V) \rightarrow (T_W,W)$ between
nonsingular toric varieties with $T_V \subset V$ and $T_W \subset W$ indicating the
tori, and points $p_V \in V$ and $q_W = \varphi(p_V) \in W$, such that we have
isomorphisms of the completions of the local rings
$$\align
\sigma:\widehat{{\Cal O}_{X,p}} &\overset{\sim}\to{\rightarrow} \widehat{{\Cal
O}_{V,p_V}} \\
\tau:\widehat{{\Cal O}_{Y,q}} &\overset{\sim}\to{\rightarrow} \widehat{{\Cal
O}_{W,q_W}} \\
\endalign$$
with a commutative diagram
$$\CD
\widehat{{\Cal O}_{X,p}} \hskip.1in @.
\underset{\sigma}\to{\overset{\sim}\to{\rightarrow}} @.
\hskip.1in \widehat{{\Cal O}_{V,p_V}} \\
@A{f^*}AA @. @A{\varphi^*}AA \\
\widehat{{\Cal O}_{Y,q}} \hskip.1in @.
\underset{\tau}\to{\overset{\sim}\to{\rightarrow}} @.
\hskip.1in \widehat{{\Cal O}_{W,q_W}} \\
\endCD$$
which induces, via inclusions, another commutative diagram among the ideals defining the
boundary divisors
$$\CD
\widehat{{\Cal I}_{D_X,p}} \hskip.1in @. \overset{\sim}\to{\rightarrow} @.
\hskip.1in \widehat{{\Cal I}_{D_V,p_V}} \\
@A{f^*}AA @. @A{\varphi^*}AA \\
\widehat{{\Cal I}_{D_Y,q}} \hskip.1in @. \overset{\sim}\to{\rightarrow} @.
\hskip.1in \widehat{{\Cal I}_{D_W,q_W}}. \\
\endCD$$

We say that a morphism $f:X \rightarrow Y$ between nonsingular varieties is
toroidal if we can introduce the structures $(U_X,X)$ and $(U_Y,Y)$ of nonsingular
toroidal embeddings, i.e., we can find divisors $D_X = X - U_X
\subset X$ and
$D_Y = Y - U_Y \subset Y$ with only simple normal crossings, such that $f:(U_X,X)
\rightarrow (U_Y,Y)$ is a toroidal morphism in the logarithmic category.

We also invite the reader to look at Abramovich-Karu[4]
Kempf-Knudsen-Mumford-Saint-Donat[21] for reference.

\vskip.1in

1.2.3. \it Definition of a morphism being log smooth. \rm \hskip.1in Let
$f:(U_X,X)
\rightarrow (U_Y,Y)$ be a morphism in the logarithmic category of nonsingular toroidal
embeddings.  We say $f$ is log smooth if the natural homomorphism
$$f^*\{\wedge^{\dim Y} \Omega_Y^1(\log D_Y)\} \wedge^{\dim X - \dim
Y}\Omega_X^1(\log D_X)\} \rightarrow \wedge^{\dim X}\Omega_X^1(\log D_X)$$ is
surjective.

\vskip.1in

1.3. $\bold{Equivalence\ of\ ``toroidal"\ and\ ``log\ smooth".}$

\vskip.1in

\proclaim{Proposition 1.3.1} Let $f:(U_X,X) \rightarrow (U_Y,Y)$ be a morphism in the
logarithmic category of nonsingular toroidal embeddings.  Then $f$ is toroidal if and
only if it is log smooth.
\endproclaim

\demo{Proof}\enddemo We only demonstrate a proof in the case where $\dim X = \dim Y =
2$.  We make a further assumption that at the point $p \in D_X \subset X$ of our
concern two irreducible components of $D_X$ meet and also that at $q = f(p)
\in D_Y
\subset Y$ two irreducible components of $D_Y$ meet.  (In the terminology of \S 3, the
points $p$ and $q$ are of type $2_p$ and $2_q$.)  Though simple, the consideration of
this typical case captures the essence of the idea.  (We refer the reader to Kato[20]
for a general proof, though it is not so much more difficult to deal with the general
case via local analysis than with the typical case.)

\vskip.1in

We choose a system of regular parameters $(x_1,x_2)$ of $\widehat{{\Cal O}_{X,p}}$ in
an analytic neighborhood
$U_p$ of
$p$ (resp. $(y_1,y_2)$ of $\widehat{{\Cal O}_{Y,q}}$ in an analytic neighborhood $U_q$
of
$q$) such that it is compatible with the logarithmic structure, i.e., $D_X \cap U_p =
\{x_1x_2 = 0\}$ (resp. $D_Y \cap U_q =
\{y_1y_2 = 0\}$).

By condition (b) $D_X = f^{-1}(D_Y)$ imposed on a morphism in the logarithmic category,
we have
$$\left\{\aligned
f^*y_1 &= u \cdot x_1^ax_2^b \\
f^*y_2 &= v \cdot x_1^cx_2^d \\
\endaligned\right.$$
where $a,b,c,d \in {\Bbb Z}_{\geq 0}$ and $u,v \in \widehat{{\Cal O}_{X,p}}^{\times}$
are units. 

We compute
$$\align
f^*(d\log y_1 \wedge d\log y_2) &= f^*\left(\frac{dy_1 \wedge dy_2}{y_1 \cdot
y_2}\right) \\
&= d\log f^*y_1 \wedge d\log f^*y_2 \\
&= \left(\frac{du}{u} + a\frac{dx_1}{x_1} + b\frac{dx_2}{x_2}\right) \wedge
\left(\frac{dv}{v} + c\frac{dx_1}{x_1} + d\frac{dx_2}{x_2}\right) \\
&= r_{\log} \cdot \frac{dx_1 \wedge dx_2}{x_1 \cdot x_2} \\
\endalign$$
where
$$r_{\log} = \det\left[\matrix
a & b \\
c & d \\
\endmatrix\right] + (\roman{higher\ terms\ \ i.e.\ =\ 0\ at\ p}).$$ 
Therefore, we conclude
$$\align
&f \roman{\ being\ log\ smooth} \\
&\Longleftrightarrow R_{\log} := \roman{div}(r_{\log})
= 0 \\
&\Longleftrightarrow \exists\ (x_1,x_2)\ \&\ (y_1,y_2) \roman{\ s.t.}
\left\{\aligned
f^*y_1 &= u \cdot x_1^ax_2^b \\
f^*y_2 &= v \cdot x_1^cx_2^d \\
\endaligned\right. \\
&\hskip.28in \roman{\ where\ }
\det\left[\matrix
a & b \\
c & d \\
\endmatrix\right] \neq 0 \roman{\ and\ }u,v \in \widehat{{\Cal O}_{X,p}}^{\times}
\roman{\ are\ units.}\\ 
& \Longleftrightarrow \exists\ (x_1,x_2)\ \&\ (y_1,y_2) \roman{\ s.t.}
\left\{\aligned
f^*y_1 &= x_1^ax_2^b \\
f^*y_2 &= x_1^cx_2^d \\
\endaligned\right. \roman{\ where\ }
\det\left[\matrix
a & b \\
c & d \\
\endmatrix\right] \neq 0 \\
&\Longleftrightarrow f \roman{\ being\ toroidal.}
\endalign$$
It is obvious that the third condition implies the second, while one can see
that the second implies the third by replacing $x_1$ with
$u^{\alpha}v^{\gamma} \cdot x_1$ and $x_2$ with $u^{\beta}v^{\delta} \cdot x_2$, where
the equation
$$\left[\matrix
a & b \\
c & d \\
\endmatrix\right]
\left[\matrix
\alpha & \gamma \\
\beta & \delta \\
\endmatrix\right] = 
\left[\matrix
1 & 0 \\
0 & 1 \\
\endmatrix\right]$$
has a (unique) solution $\left[\matrix
\alpha & \beta \\
\gamma & \delta \\
\endmatrix\right] \in M_2({\Bbb Q})$ because $\det \left[\matrix
a & b \\
c & d \\
\endmatrix\right] \neq 0$.
 
This completes the proof.

\vskip.1in

1.4. $\bold{Formulation\ of\ the\ main\ conjectures}$. 

\vskip.1in

1.4.1. \it Toroidalization Conjecture. \rm \hskip.1in Let $f:X
\rightarrow Y$ be a dominant morphism between nonsingular varieties.  Then there
should exist sequences of blowups with smooth centers $\pi_X:X' \rightarrow X$ and
$\pi_Y:Y'
\rightarrow Y$ such that the induced map $f':X' \rightarrow Y'$ is a toroidal morphism.

\vskip.1in

1.4.2. \it Resolution of singularities of morphisms in the
logarithmic category.\rm \hskip.1in Let
$f:(U_X,X)
\rightarrow (U_Y,Y)$ be a morphism in the logarithmic category of nonsingular toroidal
embeddings.  Then there exist sequences of blowups with permissible
centers $\pi_X:(U_{X'},X') \rightarrow (U_X,X)$ and
$\pi_Y:(U_{Y'},Y')
\rightarrow (U_Y,Y)$ such that the induced map $f':(U_{X'},X')  \rightarrow (U_Y',Y')$
is log smooth.  (Note that a center is called permissible if it is smooth and locally
defined by a part of a system regular parameters compatible with the log structure.)

\vskip.1in

1.5. $\bold{Conclusion}$.  An affirmative solution to the problem of resolution of
singularities of morphisms in the logarithmic category implies an affirmative solution
to the toroidalization conjecture.

In fact, let $f:X \rightarrow Y$ be a dominant morphism between nonsingular
varieties.  If $f$ is smooth, we can take $U_X = X$ and $U_Y = Y$ for the assertion of
the toroidalization conjecture.  If $f$ is not smooth, we look at the natural
homomorphism
$$f^*\{\wedge^{\dim Y}\Omega_Y^1\} \wedge^{\dim X - \dim Y}\Omega_X^1
\twoheadrightarrow {\Cal I}_R \otimes \wedge^{\dim X}\Omega_X^1 \hookrightarrow
\wedge^{\dim X}\Omega_X^1$$
where ${\Cal I}_R$ is the ideal defining the ramification locus $R$.  We take a
sequence of blowups with smooth centers $\tau_Y: \tilde{Y} \rightarrow Y$ so that
$D_{\tilde{Y}} := \tau_Y^{-1}(f(R))$ is a divisor with only simple normal crossings. 
Then we take a sequence of blowups with smooth centers $\sigma_X:\overline{X}
\rightarrow X$ so that the induced rational map
$\overline{f}:\overline{X} \rightarrow \tilde{Y}$ is a morphism.  Finally, we take a
sequence of blowups with smooth centers $\tau_X:\tilde{X} \rightarrow \overline{X}$ so
that $D_{\tilde{X}} := \tau_X^{-1}\overline{f}^{-1}D_{\tilde{Y}}$ is a divisor with only
simple normal crossings.  We have now a morphism $\tilde{f}:(U_{\tilde{X}} := \tilde{X}
\setminus D_{\tilde{X}}, \tilde{X}) \rightarrow (U_{\tilde{Y}} := \tilde{Y}
\setminus D_{\tilde{Y}}, \tilde{Y})$ in the logarithmic category of nonsingular
toroidal embeddings.  Now it is clear that a solution for resolution of singularities
of the morphism $\tilde{f}$ in the logarithmic category would imply a solution for
toroidalization of the original morphism $f$ by construction.

\newpage

\S 2. Algorithm for toroidalization

\vskip.1in

2.1. $\bold{Algorithm\ for\ toroidalization\ in\ the\ case\ where\ \dim X = \dim\ Y
= 2}$.  The purpose of this section is to describe an algorithm for toroidalization
of a dominant morphism $f:X \rightarrow Y$ between
nonsingular varieties in the case where
$\dim X = \dim Y = 2$.  We make a couple of remarks at the end of the
section regarding possible algorithms in higher dimension.

\vskip.1in

2.1.1. \it Setting.\rm \hskip.1in As in the formulation of the problem of
toroidalization in
\S 1, let
\linebreak
$f:(U_X,X)
\rightarrow (U_Y,Y)$ be a morphism
in the logarithmic category of nonsingular toroidal embeddings in the case where
$\dim X = \dim Y = 2$, that is to say: 

(a) $f:X
\rightarrow Y$ is a dominant morphism between nonsingular
varieties of

$\dim X = \dim Y = 2$, 

(b) $D_X = f^{-1}(D_Y)$, and

(c) the restriction of $f$ to the open subset $U_X$, denoted by $f|_{U_X}$, is a
smooth

morphism.

\vskip.1in

2.1.2. \it Description of the algorithm.\rm \hskip.1in Our algorithm for
toroidalization proceeds as follows:

\vskip.1in

\noindent $\boxed{\roman{Step}\ 1}$

\vskip.1in

We write down the natural homomorphism among the logarithmic differential forms
$$\align
f^*\{\wedge^{\dim Y}\Omega_Y^1(\log D_Y)\} \wedge^{\dim X - \dim
Y}\Omega_X^1(\log D_X) &\twoheadrightarrow {\Cal I}_{R_{\log}} \otimes \wedge^{\dim
X}\Omega_X^1(\log D_X) \\
&\hookrightarrow \wedge^{\dim
X}\Omega_X^1(\log D_X) \\
\endalign$$
where ${\Cal I}_{R_{\log}}$ is the ideal defining the logarithmic ramification
locus.

Note that in our case where $\dim X = \dim Y = 2$, or more generally when
\linebreak
$f:X
\rightarrow Y$ is a generically finite morphism, we have the usual logarithmic
ramification formula between the log canonical divisors
$$K_X + D_X = f^*(K_Y + D_Y) + R_{\log},$$
where $R_{\log}$ is the logarithmic ramification divisor (cf. Iitaka[17]).   The
natural homomorphism above corresponds to the inclusion of line bundles
$$\align
0 \rightarrow {\Cal O}_X(f^*(K_Y + D_Y)) &= {\Cal O}_X(- R_{\log} + K_X + D_X) \\
&\hookrightarrow {\Cal O}_X(K_X + D_X). \\
\endalign$$

\noindent $\boxed{\roman{Step}\ 2}$

\vskip.1in

We look at $f(R_{\log})$, the image of the logarithmic
ramification locus.

\vskip.1in

$\bold{The\ special\ feature,\ which\ is\ a\ consequence\ of\ the\ assumption\ \dim Y =
2,}$ 

$\bold{is\ that\ f(R_{\log})\ consists\ of\ finitely\ many\ points.}$

\vskip.1in

Blowup $Y$ with center $f(R_{\log})$ (or more precisely speaking, blowup $Y$ along
the defining ideal
${\Cal I}_{f(R_{\log})}$ of the image of the logarithmic ramification locus with the
reduced structure).

\newpage

\noindent $\boxed{\roman{Step}\ 3}$

\vskip.1in

Apply an algorithm for canonical principalization (See Remark 2.2.2 for the
precise meaning.) to the ideal
$f^{-1}({\Cal I}_{f(R_{\log})})
\cdot {\Cal O}_X$
$$\CD
(U_X,X) = (U_{X_0},X_0) @. \hskip.1in
\underset{\roman{cp}_1}\to{\overset{\roman{canonical\
principalization}}\to{\leftarrow}}
\hskip.1in @. (U_{X_1},X_1) \\ @V{f = f_0}VV  @. @V{f_1}VV \\
(U_Y,Y) = (U_{Y_0},Y_0) @. \hskip.1in
\underset{\roman{bp}_1}\to{\overset{\roman{Blowup\ }f(R_{\log})}\to{\leftarrow}}
\hskip.1in @. (U_{Y_1},Y_1) \\
\endCD$$
Note that by the universal property of blowup for principalization, the induced
rational map $f_1$ is guaranteed to be a morphism.

\vskip.1in

If the induced morphism $f_1:(U_{X_1},X_1) \rightarrow (U_{Y_1},Y_1)$ is toroidal,
we end the algorithm.

If the induced morphism $f_1:(U_{X_1},X_1) \rightarrow (U_{Y_1},Y_1)$ is not toroidal,
then we go back to $\boxed{\roman{Step}\ 1}$ of the algorithm with $f:= f_1$
and repeat the process.

\vskip.1in

2.2. $\bold{Remarks\ on\ the\ algorithm}$.

\vskip.1in

2.2.1. \it A morphism is always toroidal over $Y$ in codimension one.\rm
\hskip.1in It is straightforward to observe that, if
$\dim Y = 1$, then
$R_{\log} = \emptyset$ and hence also $f(R_{\log}) = \emptyset$. 

In fact, let $p \in D_X \subset X$ be an arbitrary closed point with $q = f(p)
\in D_Y
\subset Y$ its image.  Let $f^*:\widehat{{\Cal O}_{Y,q}} \rightarrow
\widehat{{\Cal O}_{X,p}}$ be the induced homomorphism between the completions
of the local rings.  We take some systems of regular parameters $(x_1, ... ,
x_d)$ of $\widehat{{\Cal O}_{X,p}}$ and $(y_1)$ of $\widehat{{\Cal O}_{Y,q}}$,
compatible with the logarithmic structures.  (That is to say, for some
analytic neighborhoods $U_p$ of $p$ and $U_q$ of $q$, we have
$$\align
D_X \cap U_p &= \{\prod_{m \in M} x_m = 0\}\\
D_Y \cap U_q &= \{y_1 = 0\}\\
\endalign$$
for a subset $M \subset \{1, ... , d\}$.)

\vskip.1in

We observe
$$\align
f \text{\ is\ toroidal} &\Longleftrightarrow \exists \hskip.1in (x_1, ... ,
x_d)\hskip.1in  \&
\hskip.1in (y_1) \hskip.1in \text{\ s.t.\ } \\
& \hskip.2in f^*y_1 = x_1^{a_1} \cdot\cdot\cdot x_d^{a_d} = \prod x_i^{a_i}
\text{\ with\ }a_i \geq 0 \hskip.1in \forall i \hskip.1in \text{\ and\ }
\hskip.1in a_{i_o}
\neq 0 \text{\ for\ some\ }i_o \\
&\Longleftrightarrow \exists \hskip.1in (x_1, ... , x_d)\hskip.1in  \&
\hskip.1in (y_1) \hskip.1in \text{\ s.t.\ } \\
& \hskip.2in f^*y_1 = u \cdot x_1^{a_1} ... x_d^{a_d} = u \cdot \prod x_i^{a_i} \text{\
with\ }a_i
\geq 0 \hskip.1in \forall i \hskip.1in \text{\ and\ } \hskip.1in a_{i_o} \neq 0 \text{\
for\ some\ }i_o, \\
& \hskip.2in \text{where\ }u \text{\ is\ a\ unit.}\\
\endalign$$
(It is obvious that the second condition implies the last, while one can see
that the last condition implies the second by, e.g., replacing $x_{i_o}$ with
$u^{1/i_o} \cdot x_{i_o}$, where \linebreak
$u^{1/i_o} \in \widehat{{\Cal
O}_{X,p}}$ exists because of the assumption of the base field $k$ being
algebraically closed of characteristic zero.)

\vskip.1in

The last condition obviously holds by condition (b) $D_X = f^{-1}(D_Y)$ imposed on
a morphism in the logarithmic category.

Therefore, $f$ is toroidal.

\vskip.1in

This
observation immediately implies that, for an arbitrary dominant morphism
$f:(U_X,X)
\rightarrow (U_Y,Y)$ in the logarithmic category, there is no locus of
logarithmic ramification over a generic point of codimension one, i.e., 
$$\roman{codim}_{Y}f(R_{\log}) \geq 2.$$
This is why $f(R_{\log})$ consists of finitely many points in the case where $\dim Y
= 2$.

\vskip.1in

2.2.2. \it Canonical principalization. \rm \hskip.1in By an algorithm for canonical
principalization, we mean a specific and fixed algorithm which, to a given ideal
${\Cal I}$ on a nonsingular toroidal embedding
$(U_X,X)$ with support contained in the boundary divisor, i.e., $\roman{Supp}{\Cal
O}_X/{\Cal I}
\subset D_X = X
\setminus U_X$, assigns a uniquely determined sequence of blowups
$$\align
(U_X,X) = (U_{X_0},X_0) &\overset{\pi_1}\to{\leftarrow} (U_{X_1},X_1) \leftarrow
...\\ 
... \leftarrow (U_{X_{i-1}},X_{i-1}) &\overset{\pi_i}\to{\leftarrow} (U_{X_i},X_i)
\leftarrow ... \\
... &\overset{\pi_l}\to{\leftarrow} (U_{X_l},X_l)\\
\endalign$$
with permissible centers $Y_{i-1} \subset D_{X_{i-1}} \subset X_{i-1}$ (The
adjective ``permissible" means, by definition, that $Y_{i-1}$ is smooth and that for
any closed point
$p \in Y_{i-1} \subset D_{X_{i-1}} \subset X_{i-1}$ there exists an
analytic neighborhood $U_p$ with a system of regular parameters $(x_1, ... , x_d)$
of $\widehat{{\Cal O}_{X_{i-1},p}}$ such that 
$$\align
D_{i-1} \cap U_p &= \{\prod_{m \in M}x_m = 0\} \\
Y_{i-1} \cap U_p &= \cap_{m \in N}\{x_m = 0\} \\
\endalign$$
for some subsets $M, N \subset \{1, ... , d = \dim X\}$.)

We require as a result of principalization that the
ideal $\pi^{-1}{\Cal I} \cdot {\Cal O}_{X_l}$ is principal (where $\pi = \pi_l
\circ
\cdot\cdot\cdot \circ \pi_1$), i.e., for any closed point $p \in X_l$ there
exists an analytic neighborhood $U_p$ with a system of regular
parameters $(x_1, ... , x_d)$ of $\widehat{{\Cal O}_{X_l,p}}$ such that
$$\align
D_l \cap U_p &= \{\prod_{m \in M}x_m = 0\}\\
\pi^{-1}{\Cal I} \cdot {\Cal O}_{X_l}|_{U_p} &= (\prod_{m \in L}x_m^{a_m})\\
\endalign$$
for some subsets $M, L \subset \{1, ... , d = \dim X\}$.  

\vskip.1in

The adjective ``canonical" is used to add the following two more
requirements to the algorithm: 

1. (Stability under pull-back by smooth morphisms) If
$f:(U_X,X)
\rightarrow (U_Y,Y)$ is a morphism in the logarithmic category between nonsingular
toroidal embeddings with
$f:X
\rightarrow Y$ being a smooth morphism (in the usual sense), then the process for
principalization of the ideal ${\Cal I}_X = f^*{\Cal I_Y} \subset f^*{\Cal O}_Y = {\Cal
O}_X$ specified by the algorithm on $(U_X,X)$ should coincide with the pull-back
(i.e., the one obtained by taking the Cartesian product ``$\times_YX$") of the
process for principalization of the ideal ${\Cal I}_Y$ specified by the
algorithm on
$(U_Y,Y)$.

It should be noted that this requirement of the algorithm for principalization being
stable under pull-back by smooth morphisms implies equivariance: if a group
(finite or algebraic) acts on
$(U_X,X)$ and the ideal ${\Cal I}$ is invariant under the action, then the
process of canonical principalization is equivariant, i.e., the action lifts to
the blowups in the process.  This can be seen easily by applying this requirement
to the two smooth morphisms, $pr:X
\times G
\rightarrow X$ the projectition and $\mu:X \times G \rightarrow X$ the multiplication map
representing the action, and by observing that the process of canonical principalization
on $X \times G$ coincides both with the pull-back by $pr$ and with the
pull-back by $\mu$.

2. (Analytic nature of the algorithm) Suppose that $(U_{X_{\alpha}},X_{\alpha})$
with
${\Cal I}_{\alpha}$ and $(U_{X_{\beta}},X_{\beta})$ with ${\Cal
I}_{\beta}$ are gievn as above.  Suppose that we have
closed points $p_{\alpha}
\in X_{\alpha}$ and $p_{\beta} \in X_{\beta}$
with an isomorphism between the completions 
$$\widehat{{\Cal O}_{X_{\alpha},p_{\alpha}}} = {\Cal
O}_{X_{\alpha},p_{\alpha}} \otimes \widehat{{\Cal O}_{X_1,p_1}}
\overset{\sim}\to{\rightarrow} {\Cal O}_{X_{\beta},p_{\beta}} \otimes
\widehat{{\Cal O}_{X_{\beta},p_{\beta}}} =
\widehat{{\Cal O}_{X_{\beta},p_{\beta}}}$$ 
which induces isomorphisms via inclusions
$$\align
\widehat{{\Cal I}_{\alpha}} = {\Cal I}_{\alpha} \otimes \widehat{{\Cal
O}_{X_{\alpha},p_{\alpha}}} &
\overset{\sim}\to{\rightarrow} {\Cal I}_{\beta} \otimes \widehat{{\Cal
O}_{X_{\beta},p_{\beta}}} = \widehat{{\Cal I}_{\beta}}\\
\widehat{{\Cal I}_{D_{\alpha}}} = {\Cal I}_{D_{\alpha}} \otimes \widehat{{\Cal
O}_{X_{\alpha},p_{\alpha}}} &
\overset{\sim}\to{\rightarrow} {\Cal I}_{D_{\beta}} \otimes \widehat{{\Cal
O}_{X_{\beta},p_{\beta}}} = \widehat{{\Cal I}_{D_{\beta}}}.\\
\endalign$$
Then the isomorphism
$$(U_{X_{\alpha}},X_{\alpha}) \times \roman{Spec}\ \widehat{{\Cal
O}_{X_{\alpha},p_{\alpha}}}
\overset{\sim}\to{\rightarrow} (U_{X_{\beta}},X_{\beta}) \times \roman{Spec}\
\widehat{{\Cal O}_{X_{\beta},p_{\beta}}}$$
lifts to the isomorphisms between the blowups in the process of principalization
for the ideal
${\Cal I}_{\alpha}$ on $(U_{X_{\alpha}},X_{\alpha})$, pulled back by $\roman{Spec}\
\widehat{{\Cal O}_{X_{\alpha},p_{\alpha}}}$, and the blowups in the process of
principalization for the ideal
${\Cal I}_{\beta}$ on $(U_{X_{\beta}},X_{\beta})$, pulled back by $\roman{Spec}\
\widehat{{\Cal O}_{X_{\beta},p_{\beta}}}$.

3. (Avoidance of blowing up the points where the ideal is already principal) If
the ideal ${\Cal I}$ is already principal at a closed point $p \in X$, then the
algorithm for canonical principalization leaves a neighborhood of $p$ untouched. 

\vskip.1in

An important consequence of requirements 1 and 2 is that if we
apply the algorithm for canonical principalization to a toroidal ideal, then the
blowups in the process of principalization are all toroidal as well as the
ideals which are the total transforms of the original toroidal ideal.

\vskip.1in

\it In our special case where $\dim X = \dim Y = 2$, there is a unique algorithm for
canonical principalization, namely, at each stage on $X_{i-1}$ we only blow up the
points where the total transform ${\Cal I}_{i-1}$ of the ideal is \it not \rm
principal.\rm

\vskip.1in

2.2.3 \it Notation for the logarithmic ramification locus.\rm \hskip.1in Given a
morphism $f:(U_X,X) \rightarrow (U_Y,Y)$ in the logarithmic category, if we want
to emphasize that the logarithmic ramification locus
$R_{\log}$ is associated to the morphism $f$, we write $R_{\log,f}$ instead of
$R_{\log}$.

\newpage  

2.3. $\bold{Further\ remarks\ on\ the\ algorithm}$

\vskip.1in

2.3.1. \it Case where $\dim X > 2$ and $\dim Y = 2$ \rm \hskip.1in Our algorithm
formally makes sense even in the case $\dim X > 2$ as long as $\dim Y = 2$,
adopting any algorithm for canonical principalization, such as the one prescribed
by Bierstone-Milman[8] or Encinas-Villamayor[14][15] among others, satisfying the
above requirements.  However, the author cannot show that the algorithm terminates
even in the case $\dim X = 3$.

The major difference between the case where $\dim X > 2$ and the case where $\dim
X = 2$ is that in the latter the defining ideal for the logarithmic ramification
locus is always principal whereas in the former may not be.  In fact, even when
$\dim X > 2$, if the defining ideal for the logarithmic ramification
locus happens to be principal,  a combinatorial analysis very similar to the one in
the case where $\dim X = 2$ would yield toroidalization.  Therefore, one could
try to formulate an algorithm where we first aim at principalizing the defining
ideal for the logarithmic ramification locus.  An interpretation of the monumental
work by Cutkosky[12] could be given from this point of view.  Deatails will be
published elsewhere.

\vskip.1in

2.3.2. \it Case where $\dim Y > 2$ \rm \hskip.1in We are completely at loss in the
case where $\dim Y > 2$, even in the case where $\dim X = \dim Y = 3$ and hence the
defining ideal for the logarithmic ramification locus is principal.  The stalemate
is a result of our (or rather, the author's) lack of understanding of the expected
inductional structure in the scheme of toroidalization, though the logarithmic
category should be tailor-made for such an inductional scheme via adjunction.

\newpage

\S 3. Local analysis in the case where $\dim X = \dim Y = 2$

\vskip.1in

3.1. $\bold{Purpose\ of\ this\ section}$.  The purpose of this section is to
present a detailed local analysis of what happens in $\boxed{\roman{Steps}\
1, 2, 3}$ of the algorithm 2.1 for toroidalization in the case where $\dim X
= \dim Y = 2$, under the same setting as in 2.1.1, i.e., $f:(U_X,X)
\rightarrow (U_Y,Y)$ is a morphism in the logarithimic category of
nonsingular toroidal embeddings in the case where $\dim X = \dim Y = 2$ so
that

(a) $f:X
\rightarrow Y$ is a dominant morphism between nonsingular
varieties of

$\dim X = \dim Y = 2$, 

(b) $D_X = f^{-1}(D_Y)$, and

(c) the restriction of $f$ to the open subset $U_X$, denoted by $f|_{U_X}$, is a
smooth

morphism.

\vskip.1in

3.2. \it Type and Case Descriptions.\rm \hskip.1in Let $p \in D_X \subset X$
be a closed point and $q = f(p) \in D_Y \subset Y$ its image.  In the
following analysis, we always choose $q$ to be in the image of the
logarithmic ramification locus, i.e., $q \in f(R_{\log})$.

\vskip.1in

We make our analysis
sorting out the types and cases, based upon the following criteria:

\vskip.1in

$\circ$ The morphism $f$ is toroidal or not at $p$.  

\hskip.2in We put the letter $T$ when
it is toroidal, and the letter $N$ when not.

$\circ$ Description of the boundary divisors at $p$ and $q$ and their
behavior with respect 

to $f$.

\vskip.1in

We have the following type descriptions on the boundary divisors at $p$ and
$q$ (We take some systems of regular parameters $(x_1,x_2)$ of
$\widehat{{\Cal O}_{X,p}}$ and $(y_1,y_2)$ of $\widehat{{\Cal O}_{Y,q}}$,
compatible with the logarithmic structures, in analytic neighborhoods
$U_p$ of $p$ and $U_q$ of $q$, respectively.):

\vskip.1in

$2_p$: The point $p$ is the intersection point of two irreducible components

\hskip.3in $G_1 =
\{x_1 = 0\}$ and $G_2 = \{x_2 = 0\}$ with $G_1 \cup G_2 = D_X
\cap U_p$ 

\hskip.3in for some analytic neighborhood $U_p$ of $p$.

$1_p$: The point $p$ belongs to only one irreducible component $G_1 =
\{x_1 = 0\}$ 

\hskip.3in with $G_1 = D_X \cap U_p$ for some analytic neighborhood $U_p$ of
$p$, while 

\hskip.3in $G_2 = \{x_2 = 0\} \not\subset D_X \cap U_p$.

$2_q$: The point $q$ is the intersection point of two irreducible components

\hskip.3in $H_1 = \{y_1 = 0\}$ and $H_2 = \{y_2 = 0\}$ with $H_1 \cup H_2 = D_Y
\cap U_q$ 

\hskip.3in for some analytic neighborhood $U_q$ of $q$

$1_q$: The point $q$ belongs to only one irreducible component $H_1 =
\{y_1 = 0\}$ 

\hskip.3in with $H_1 = D_Y \cap U_q$ for some analytic neighborhood $U_q$ of
$q$, while 

\hskip.3in $H_2 = \{y_2 = 0\} \not\subset D_Y \cap U_q$.

\newpage

We have the following ten subcases (some of which, namely, Subcase
$2_p1_q0$ and Subcase $1_p2_q0$, will be doomed impossible), where the third number
indicates how many irreducible components of the boundary divisor map onto $q$: 

\vskip.1in

\hskip.1in Subcase $2_p2_q0$: Neither $G_1$ nor $G_2$ maps onto
$q$

\hskip.1in Subcase $2_p2_q1$: Only one of $G_1$ or $G_2$ (say, $G_1$) maps onto
$q$

\hskip.1in Subcase $2_p2_q2$: Both $G_1$ and $G_2$ maps onto $q$

\hskip.1in Subcase $1_p2_q0$: $G_1$ does not map onto $q$ (and hence maps onto
$H_1$ (or $H_2$))

\hskip.1in Subcase $1_p2_q1$: $G_1$ maps onto $q$

\hskip.1in Subcase $1_p1_q0$: $G_1$ does not map onto $q$ (and hence $G_1$ maps
onto $H_1$)

\hskip.1in Subcase $1_p1_q1$: $G_1$ maps onto $q$

\hskip.1in Subcase $2_p1_q0$: Neither $G_1$ nor $G_2$ maps onto
$q$

\hskip.1in Subcase $2_p1_q1$: Only one of $G_1$ or $G_2$ (say, $G_1$) maps onto
$q$

\hskip.1in Subcase $2_p1_q2$: Both $G_1$ and $G_2$ maps onto $q$

\vskip.1in

3.3. \it Case by case analysis.\rm \hskip.1in Now we carry out a detailed
analysis, in each (sub)case, of what happens when we blowup $q$ in
$\boxed{\roman{Step}\ 2}$ and when we take the canonical principalization of
$f^{-1}{\Cal I}_q \cdot {\Cal O}_X = f^{-1}m_q
\cdot {\Cal O}_X$ in $\boxed{\roman{Step}\ 3}$.

\vskip.1in

We provide:

\vskip.1in

A) Coordinate Expression

B) Descriptions of the canonical principalization $\roman{cp}_1$ and the induced

\hskip.2in morphism $f_1$

C) Conclusion on the behavior of the logarithmic ramification divisor

\vskip.2in

It turns out that what weighs more, for the purpose of sorting out the
subcases, is the type description of the point
$q$ than that of $p$, and we regroup the subcases accordingly.

\vskip.1in

Remark on notation: When we write Subcases $* \hskip.03in 2_q \hskip.03in
*
\hskip.03in T$, we indicate that we are dealing with all the subcases of
type $2_p2_q0$,
$2_p2_q1$, $2_p2_q2$, $1_p2_q0$, $1_p2_q1$, collectively, with the extra
condition that $f$ is toroidal at $p$. 

When we write $1_p2_q0$ without mentioning whether $f$ is toroidal or not at
$p$, we indicate that we are dealing with all the possibilities in the
subcase without a priori putting the extra assumption on whether $f$ is
toroidal or not.

\newpage

3.3.1) Subcases $* \hskip.03in 2_q \hskip.03in *$: These consist of the
following subcases:

\vskip.1in

\hskip.1in Subcase $2_p2_q0$

\hskip.1in Subcase $2_p2_q1$

\hskip.1in Subcase $2_p2_q2$

\hskip.1in Subcase $1_p2_q1$

\vskip.1in

\hskip.1in Subcase $1_p2_q0$: This subcase does not happen.  In fact, say $G_1$
maps onto $H_1$ (and hence does not map onto $H_2$).  Note that
in a neighborhood
$U_p$ of $p$ we have $f^{-1}(H_2) \cap U_p
\subset f^{-1}(D_Y) \cap U_p = D_X \cap U_p = G_1$.  Now since
$f^{-1}(H_2)$ is of pure codimension one, being the pull-back of a Cartier
divisor, and since $p \in f^{-1}(H_2)$, we conclude $G_1 \cap U_p =
f^{-1}(H_2) \cap U_p$.  But this would imply either $G_1$ maps onto $q$ or
maps onto
$H_2$, a contradiction!

\vskip.1in

3.3.1.A) Coordinate Expression

\vskip.1in

By the computation presented in the proof of Proposition 1.3.1, we conclude
the following.

\vskip.1in

Subcases $2_p2_q \hskip.03in *T$: There exist systems of regular parameters
$(x_1,x_2)$ of $\widehat{{\Cal O}_{X,p}}$ and $(y_1,y_2)$ of $\widehat{{\Cal
O}_{Y,q}}$ compatible with the logarithmic structures such that
$$\left\{\aligned
f^*y_1 &= x_1^ax_2^b \\
f^*y_2 &= x_1^cx_2^d \\
\endaligned
\right. \hskip.1in \text{where} \hskip.1in \det\left[\matrix
a & b \\
c & d \\
\endmatrix\right] \neq 0.$$

Subcases $2_p2_q \hskip.03in *N$: There exist systems of regular parameters
$(x_1,x_2)$ of $\widehat{{\Cal O}_{X,p}}$ and $(y_1,y_2)$ of $\widehat{{\Cal
O}_{Y,q}}$ compatible with the logarithmic structures such that
$$\left\{\aligned
f^*y_1 &= u \cdot x_1^ax_2^b \\
f^*y_2 &= v \cdot x_1^cx_2^d \\
\endaligned
\right. \hskip.1in \text{where} \hskip.1in \det\left[\matrix
a & b \\
c & d \\
\endmatrix\right] = 0 \hskip.1in \text{with\ units} \hskip.1in u,v \in
\widehat{{\Cal O}_{X,p}}^{\times}.$$

Subcases $1_p2_q \hskip.03in *$: There exist systems of regular parameters
$(x_1,x_2)$ of $\widehat{{\Cal O}_{X,p}}$ and $(y_1,y_2)$ of $\widehat{{\Cal
O}_{Y,q}}$ compatible with the logarithmic structures, as described in
3.2, such that
$$\left\{\aligned
f^*y_1 &= u \cdot x_1^a \\
f^*y_2 &= v \cdot x_1^c \\
\endaligned
\right. \hskip.1in \text{where\ } a,c \in {\Bbb Z}_{> 0}
\text{\ with\ units} \hskip.1in u,v \in
\widehat{{\Cal O}_{X,p}}^{\times}.$$

\vskip.1in

3.3.1.B) Descriptions of the canonical principalization $\roman{cp}_1$ and
the induced 

\hskip.3in morphism $f_1$

\vskip.1in

Subcases $2_p2_q \hskip.03in *T$: We observe that the ideal ${\Cal I}_q = m_q
= (y_1,y_2)$ is toroidal \linebreak
(in a neighborhood of $q$) in these
subcases and that the ideal \linebreak
$f^{-1}(y_1,y_2) \cdot {\Cal O}_{X,p} =
(x_1^ax_2^b,x_1^cx_2^d)$ is also toroidal (in a neighborhood of $p$). 
Therefore, by remark 2.2.2, the blowup
$\roman{bp}_1:(U_{Y_1},Y_1) \rightarrow (U_{Y_0},Y_0)$ is toroidal over a neighborhood
of $q$, and so is the canonical principalization
$\roman{cp}_1:(U_{X_1},X_1)
\rightarrow (U_{X_0},X_0) = (U_X,X)$ over a neighborhood of $p$.  Since a
composite of toroidal morphisms is again toroidal and
so is a composite of a toroidal morphism with the inverse of a toroidal
birational morphism (assuming that the composite is again a morphism, which
is the case here), we conclude that
$f_1 = \roman{bp}_1^{-1} \circ f_0 \circ \roman{cp}_1$ stays toroidal on the
locus of
$(U_{X_1},X_1)$ which sits over a neighborhood of $p$.

\vskip.1in

Subcases $2_p2_q \hskip.03in *N$: In these subcases, as we observed, we have
$$\det\left[\matrix
a & b \\
c & d \\
\endmatrix\right] = 0,$$
that is to say, the vectors $(a,b), (c,d) \in {\Bbb Z}_{\geq 0} \times {\Bbb Z}_{\geq
0}$ are linearly dependent.  (We may choose $(y_1,y_2)$, without loss of generality,
so that
$a
\geq c, b \geq d$.  Note that
$(a,b)
\neq (0,0)$ and
$(c,d)
\neq
 (0,0)$.)  This implies that the ideal $f^{-1}(y_1,y_2)
\cdot {\Cal O}_{X,p} = (x_1^cx_2^d)$
is already principal in
a neighborhood of $p$.  Therefore, by remark 2.2.2 the canonical
principalization
$\roman{cp}_1$ leaves a neighborhood of $p$ untouched.  Also, by remark
2.2.2, in a neighborhood of $p$ the map $f_1 = \roman{bp}_1^{-1} \circ f_0
\circ
\roman{cp}_1 = \roman{bp}_1^{-1} \circ f_0$ is a morphism.

\vskip.1in

Subcases $1_p2_q \hskip.03in *$: In these subcases, the ideal
$f^{-1}(y_1,y_2)
\cdot {\Cal O}_{X,p} = (x_1^a)$ (We may choose $(y_1,y_2)$ without loss of
generality so that $a
\leq c$ in the coordinate expression.) is already
principal in a neighborhood of $p$.  Therefore, by remark 2.2.2 the canonical
principalization
$\roman{cp}_1$ leaves a neighborhood of $p$ untouched.  Also, by remark
2.2.2, in a neighborhood of $p$ the map $f_1 = \roman{bp}_1^{-1} \circ f_0
\circ
\roman{cp}_1 = \roman{bp}_1^{-1} \circ f_0$ is a morphism.

\vskip.1in

3.3.1.C) Conclusion on the behavior of the logarithmic ramification divisor

\vskip.1in

In each of the subcases above, we conclude that over a neighborhood of the
point $q$, since $\roman{bp}_1$ is toroidal,
$$K_{Y_1} + D_{Y_1} = \roman{bp}_1^*(K_{Y_0} + D_{Y_0}) \hskip.1in
\text{i.e.,}
\hskip.1in \wedge^{\dim Y_1} \Omega_{Y_1}^1(\log D_{Y_1}) = \roman{bp}_1^*\wedge^{\dim Y_0}
\Omega_{Y_0}^1(\log D_{Y_0}),$$ 
and hence that the
logarithmic ramification divisor stays intact, i.e.,
$$R_{\log,f_0} = R_{\log,f_1}.$$

\newpage

3.3.2) Subcase  $1_p1_q0$: $G_1$ does not map onto $q$ (and hence $G_1$ maps
onto $H_1$)

\vskip.1in

3.3.2.A) Coordinate Expression

\vskip.1in

In this subcase, we show that $f$ is necessarily toroidal at
$p$.  (Although this is essentially due to Abhyankar's lemma, we will not use
it explicitly).  We will show by some elementary calculation that there
exist systems of regular parameters
$(x_1,x_2)$ of
$\widehat{{\Cal O}_{X,p}}$ and
$(y_1,y_2)$ on
$\widehat{{\Cal O}_{Y,q}}$, compatible with the logarithmic structures, such
that 
$$\left\{\aligned
f^*y_1 &= x_1^a \\
f^*y_2 &= x_2\\
\endaligned\right. \hskip.1in \text{with} \hskip.1in a > 0$$
using condition (c) $f|_{U_X}$ being smooth, imposed on a morphism in the
logarithmic category (cf. 1.2.1).

First, we start with some systems of regular parameters as described in
3.2.  Then since $f^{-1}(D_Y) = D_X$ where $D_Y \cap U_q = H_1 = \{y_1 = 0\}$
and
$D_X
\cap U_p = G_1 = \{x_1 = 0\}$, we have
$$f^*y_1 = u \cdot x_1^a \hskip.1in \text{for\ some\ unit\ }u \in
\widehat{{\Cal O}_{X,p}}^{\times}\hskip.1in
\text{with}
\hskip.1in a > 0.$$
By replacing $x_1$ with $u^{1/a} \cdot x_1$, we may assume
$$f^*y_1 = x_1^a.$$ 
Set
$$f^*y_2 = \Sigma_{i > 0, j \geq 0, \det\left[\matrix
a & 0 \\
i & j \\
\endmatrix\right] = 0}\alpha_{ij}x_1^ix_2^j + \Sigma_{i \geq 0, j > 0,
\det\left[\matrix a & 0 \\
i & j \\
\endmatrix\right] \neq 0}\alpha_{ij}x_1^ix_2^j,$$
where there exists $j > 0$ with $\alpha_{0j} \neq 0$, since $x_1$ should not divide
$f^*y_2$ as $G_1 = \{x_1 = 0\}$ does not map onto $q$ by the subcase
assumption.

We set
$$j_o := \min\{j;\alpha_{0j} \neq 0\}.$$

The basic point of computing the ramification is the simple observation that
$$d(x_1^ax_2^b) \wedge d(x_1^cx_2^d) = \det\left[\matrix
a & b \\
c & d \\
\endmatrix\right] \cdot x_1^{a + c - 1}x_2^{b + d -1} dx_1 \wedge dx_2.$$
We compute
$$\align
f^*(dy_1 \wedge dy_2) &= df^*y_1 \wedge df^*y_2 \\
&= d(x_1^a) \wedge d\left(\Sigma_{i > 0, j \geq 0, \det\left[\matrix
a & 0 \\
i & j \\
\endmatrix\right] = 0}\alpha_{ij}x_1^ix_2^j + \Sigma_{i \geq 0, j > 0,
\det\left[\matrix a & 0 \\
i & j \\
\endmatrix\right] \neq 0}\alpha_{ij}x_1^ix_2^j\right)\\
&= d(x_1^a) \wedge d\left(\Sigma_{i \geq 0, j > 0, \det\left[\matrix
a & 0 \\
i & j \\
\endmatrix\right] \neq 0}\alpha_{ij}x_1^ix_2^j\right)\\
&= (\Sigma_{i \geq 0, j > 0, \det\left[\matrix
a & 0 \\
i & j \\
\endmatrix\right] \neq 0} \det\left[\matrix
a & 0 \\
i & j \\
\endmatrix\right] \cdot \alpha_{ij}x_1^{a + i - 1}x_2^{j - 1}) \cdot dx_1 \wedge dx_2
\\ &= x_1^{a - 1} \cdot \left(\Sigma_{i \geq 0, j > 0, \det\left[\matrix
a & 0 \\
i & j \\
\endmatrix\right] \neq 0} \det\left[\matrix
a & 0 \\
i & j \\
\endmatrix\right] \cdot \alpha_{ij}x_1^ix_2^{j - 1}\right) \cdot dx_1 \wedge dx_2.
\endalign$$
Therefore, if $j_o > 1$, then $f$ ramifies along 
$$\{\left(\Sigma_{i \geq 0, j > 0, \det\left[\matrix
a & 0 \\
i & j \\
\endmatrix\right] \neq 0} \det\left[\matrix
a & 0 \\
i & j \\
\endmatrix\right] \cdot \alpha_{ij}x_1^ix_2^{j - 1}\right) = 0\}$$ 
other than along $\{x_1 = 0\}$,
a contradiction to condition $f|_{U_X}$ being smooth!  Therefore, we
conclude
$$j_o = 1.$$
By replacing $x_2$ with $f^*y_2$, we obtain the desired systems of regular
parameters.

Now since neither the divisor $\{y_2 = 0\}$ nor $\{x_2 = 0\}$ belongs to the
boundary divisor defining the logarithmic structure, it is clear $f$ is
toroidal at $p$ from the coordinate expression.

\vskip.1in

3.3.2.B) Descriptions of the canonical principalization $\roman{cp}_1$ and
the induced 

\hskip.3in morphism $f_1$

\vskip.1in

We consider the 1st blowup $\roman{cp}_{1,1}$ of the canonical principalization
$\roman{cp}_1$, which is a sequence $\roman{cp}_1 = \roman{cp}_{1,1} \circ
\cdot\cdot\cdot \circ \roman{cp}_{1,l}$ of blowups specified by the canonical
principalization algorithm.

\proclaim{Diagram 3.3.2.B.1}\endproclaim

\vskip2.5in

\proclaim{Claim 3.3.2.B.2} 

a) The rational map $f_{1,1} = \roman{bp}^{-1} \circ f_0
\circ
\roman{cp}_{1,1}$ is well-defined (regular), that is to say, the ideal
$(f_0 \circ \roman{cp}_{1,1})^{-1}(y_1,y_2) \cdot {\Cal O}_{X_{1,1}}$ is
principal, except possibly at $p_1'$.

b) At $p_1'$, the morphism $f \circ \roman{cp}_{1,1}$ is in Subcase
$1_{p_1'}1_q1$.
\endproclaim

\demo{Proof}\enddemo a) Observe that the morphisms $f_0$ and
$\roman{cp}_{1,1}$ and the ideal $(y_1,y_2)$ are toroidal, with respect to the
modified logarithmic structures obtained by adding
$\{y_2 = 0\}$ and
$\{x_2 = 0\}$ to the original boundary divisors.  From this it follows easily
that $f_{1,1}$ is regular, except possibly at $p_1$ and/or $p_1'$.

Now at $p_1$, we have a system of regular parameters $(\frac{x_1}{x_2},x_2)$
with coordinate expression
$$\left\{\aligned
(f_0 \circ \roman{cp}_{1,1})^*y_1 &= ((\frac{x_1}{x_2})x_2)^a \\
(f_0 \circ \roman{cp}_{1,1})^*y_2 &= x_2, \\
\endaligned\right.$$
which immediately implies that the ideal $(f_0 \circ
\roman{cp}_{1,1})^{-1}(y_1,y_2) \cdot \widehat{{\Cal O}_{X_{1,1},p_1}} =
(x_2)$ is principal at
$p_1$ and hence that
$f_{1,1}$ is regular at $p_1$.

b) The verification for statement b) is immediate.

\vskip.1in

3.3.2.C) Conclusion on the behavior of the logarithmic ramification divisor

\proclaim{Claim 3.3.2.C.1}

a) The coefficient of $G_1$ in the logarithmic ramification divisor remains 0, i.e.,
$$\nu_{G_1}(R_{\log,f_0}) = \nu_{G_1}(R_{\log,f_1}) = 0.$$
b) The strict transform of the exceptional divisor $E_p$ for $\roman{cp}_{1,1}$,
obtained by blowing up $p$, does not appear in $R_{\log,f_1}$, i.e.,
$$\nu_{E_p}(R_{\log,f_1}) = 0.$$
More generally, none of the irreducible components $E$ of the exceptional
divisor for
$\roman{cp}_1$ (i.e., none of the strict transforms of the exceptional
divisors for
$\roman{cp}_{1,1}, ... ,\roman{cp}_{1,l}$) appear in
$R_{\log,f_1}$, i.e.,
$$\nu_E(R_{\log,f_1}) = 0.$$
\endproclaim

\demo{Proof}\enddemo a) The verification for statement a) is obvious.

b) Take the standard neighborhood of $q_1$ with a system
of regular parameters $(\frac{y_1}{y_2},y_2)$.  Then observing
$$\left\{\aligned
\nu_{E_p}(\frac{y_1}{y_2}) &= a - 1 \geq 0 \\
\nu_{E_p}(y_2) &= 1 > 0,\\
\endaligned\right.$$
we conclude that the generic point of $E_p$ maps into the standard neighborhood of
$q_1$.

Observe that the morphisms $f_0, \roman{bp}_1, \roman{cp}_{1,1}$ are all
toroidal with respect to the modified logarithmic structures obtained by
adding
$\{y_2 = 0\}$ and $\{x_2  = 0\}$ (and their pull-backs) to the original boundary
divisors.  The original logarithmic structures coincide with the modified ones
in a neighborhood of the generic point of $E_p$ and in the standard neighborhood
of $q_1$.  Therefore, $f_{1,1}$ is toroidal in a neighborhood of the generic
point of $E_p$ with respect to the original logarithmic structures, and hence
$$\nu_{E_p}(R_{\log,f_1}) = 0.$$
This proves the first part of statement b).

The second part of statement b) follows from Claim 3.3.2.B.2 b) and
\linebreak Claim
3.3.3.C.1 b).

\newpage

3.3.3) Subcase $1_p1_q1$: $G_1$ maps onto $q$

\vskip.1in

3.3.3.A) Coordinate Expression

\vskip.1in

First, we start with some systems of regular parameters as described in
3.2.  Then since $f^{-1}(D_Y) = D_X$ where $D_Y \cap U_q = H_1 = \{y_1 = 0\}$
and
$D_X
\cap U_p = G_1 = \{x_1 = 0\}$, we have
$$f^*y_1 = u \cdot x_1^a \hskip.1in \text{for\ some\ unit\ }u \in
\widehat{{\Cal O}_{X,p}}^{\times}\hskip.1in
\text{with}
\hskip.1in a > 0.$$
By replacing $x_1$ with $u^{1/a} \cdot x_1$, we may assume
$$f^*y_1 = x_1^a.$$
Set
$$f^*y_2 = \Sigma_{i > 0, j \geq 0, \det\left[\matrix
a & 0 \\
i & j \\
\endmatrix\right] = 0}\alpha_{ij}x_1^ix_2^j + \Sigma_{i \geq 0,
j > 0, \det\left[\matrix
a & 0 \\
i & j \\
\endmatrix\right] \neq 0}\alpha_{ij}x_1^ix_2^j,$$
where no term of the form $x_2^j = x_1^0x_2^j\hskip.1in (i = 0, j > 0)$ appears,
since
$x_1$ has to divide $f^*y_2$.

Now since $f|_{U_X}$ must be smooth by condition (c) imposed on a morphism
in the logarithmic category (cf. 1.2.1) and by computing the
ramification \linebreak
$f^*(dy_1
\wedge dy_2)/dx_1 \wedge dx_2$ as in 3.3.2.A, we conclude that $f^*y_2$
must be of the form
$${\aligned
f^*y_2 &= \Sigma_{i \geq i_s > 0, j = 0}\alpha_{i0}x_1^i \\
&+
\alpha_{i_o,1}x_1^{i_o}x_2^1 + \Sigma_{i \geq i_o > 0,j >
0, (i,j) \neq (i_o,1)}\alpha_{ij}x_1^ix_2^j\\
\endaligned}  \hskip.1in \text{with} \hskip.1in \alpha_{i_o,1} \neq 0,$$
where
$$\align
i_o &= \min\{i;\alpha_{ij} \neq 0, \det\left[\matrix
a & 0 \\
i & j \\
\endmatrix\right] \neq 0\} \\
i_s &= \min\{i;\alpha_{ij} \neq 0, \det\left[\matrix
a & 0 \\
i & j \\
\endmatrix\right] = 0\}.\\
\endalign$$

3.3.3.B) Descriptions of the canonical principalization $\roman{cp}_1$ and
the induced 

\hskip.3in morphism $f_1$

\vskip.1in

Case $i_s \leq i_o$: In this case, we have
$$\left\{\aligned
f^*y_1 &= x_1^a \\
f^*y_2 &= v \cdot x_1^{i_s} \\
\endaligned\right. \hskip.1in \text{with} \hskip.1in v \in
\widehat{{\Cal O}_{X,p}}^{\times}\hskip.1in
\text{unit},$$ which implies $f^{-1}(y_1,y_2) \cdot \widehat{{\Cal O}_{X,p}}
= (x_1^a, x_1^{i_s}) = (x_1^{\min\{a,i_s\}})$ is principal.  Therefore,
$\roman{cp}_1$ is an isomorphism in a neighborhood of $p$.

\vskip.1in

Case $i_s > i_o$: In this case, we have
$$\left\{\aligned
f^*y_1 &= x_1^a \\
f^*y_2 &= x_1^{i_o} \cdot (\Sigma_{i \geq i_s > 0, j = 0}\alpha_{i0}x_1^{i - i_o} +
\alpha_{i_o,1}x_2 + \Sigma_{i \geq i_o, j > 0, (i,j) \neq
(i_o,1)}\alpha_{ij}x_1^{i-i_o}x_2^j). \\
\endaligned\right.$$
Replacing $x_2$ by $(\Sigma_{i \geq i_s > 0, j = 0}\alpha_{i0}x_1^{i - i_o} +
\alpha_{i_o,1}x_2 + \Sigma_{i \geq i_o, j > 0, (i,j) \neq
(i_o,1)}\alpha_{ij}x_1^{i-i_o}x_2^j)$, we have
$$\left\{\aligned
f^*y_1 &= x_1^a \\
f^*y_2 &= x_1^cx_2 \\
\endaligned\right. \hskip.1in \text{with} \hskip.1in c = i_o.$$

\hskip.1in Subcase $a < c+1$, i.e., $a \leq c$: In this subcase, the ideal
$f^{-1}(y_1,y_2) \cdot \widehat{{\Cal O}_{X,p}} = (x_1^a)$ is principal. 
Therefore,
$\roman{cp}_1$ is an isomorphism in a neighborhood of
$p$.

\vskip.1in

\hskip.1in Subcase $a \geq c + 1$: We consider the 1st blowup
$\roman{cp}_{1,1}$ of the canonical principalization
$\roman{cp}_1$, which is a sequence of blowups $\roman{cp}_1 =
\roman{cp}_{1,1} \circ \cdot\cdot\cdot \circ \roman{cp}_{1,l}$ specified by the
canonical principalization algorithm.

\proclaim{Diagram 3.3.3.B.1}\endproclaim

\vskip2.5in

\proclaim{Claim 3.3.3.B.2} The canonical principalization $\roman{cp}_1$ is
an isomorphism in a neighborhood of $p$, except under the subcase $a \geq c +
1$ in the case
$i_o > i_1$.  In the exceptional subcase, we claim:

a) The rational map $f_{1,1} = \roman{bp}^{-1} \circ f_0
\circ
\roman{cp}_{1,1}$ is well-defined (regular), that is to say, the ideal
$(f_0 \circ \roman{cp}_{1,1})^{-1}(y_1,y_2) \cdot {\Cal O}_{X_{1,1}}$ is
principal, except possibly at $p_1'$.

b) At $p_1'$, the morphism $f \circ \roman{cp}_{1,1}$ is in Subcase
$1_{p_1'}1_q1$.
\endproclaim

\demo{Proof}\enddemo The first part of the claim is already verified.  So we have only
to prove the assertions under the subcase $a \geq c + 1$ in the case $i_o > i_1$.

a) Observe that the morphisms $f_0$ and
$\roman{cp}_{1,1}$ and the ideal $(y_1,y_2)$ are toroidal, with respect to the
modified logarithmic structures obtained by adding
$\{y_2 = 0\}$ and
$\{x_2 = 0\}$ to the original boundary divisors.  From this it follows easily
that $f_{1,1}$ is regular, except possibly at $p_1$ and/or $p_1'$.

Now at $p_1$, we have a system of regular parameters $(\frac{x_1}{x_2},x_2)$
with coordinate expression
$$\left\{\aligned
(f_0 \circ \roman{cp}_{1,1})^*y_1 &= ((\frac{x_1}{x_2})x_2)^a =
(\frac{x_1}{x_2})^ax_2^a\\ 
(f_0 \circ \roman{cp}_{1,1})^*y_2 &= ((\frac{x_1}{x_2})x_2)^cx_2 =
(\frac{x_1}{x_2})^cx_2^{c+1}\\
\endaligned\right.$$
which immediately implies that the ideal $(f_0 \circ
\roman{cp}_{1,1})^{-1}(y_1,y_2) = ((\frac{x_1}{x_2})^cx_2^{c+1})$ is principal
at
$p_1$ and hence that
$f_{1,1}$ is regular at $p_1$.

b) The verification for statement b) is immediate.

\vskip.1in

3.3.3.C) Conclusion on the behavior of the logarithmic ramification divisor

\vskip.1in

\proclaim{Claim 3.3.3.C.1} a) The coeffcient of $G_1$ in the logarithmic
ramification divisor strictly drops.  More precisely, 
$$\nu_{G_1}(R_{\log,f_1}) = \left\{\aligned
& i_o - \min\{i_o,i_s\} \hskip.1in \text{if} \hskip.1in a \geq \min\{i_o,i_s\}
\\ 
& i_o - a \hskip.67in \text{if} \hskip.1in a < \min\{i_o,i_s\}
\\ 
\endaligned\right\} < i_o = \nu_{G_1}(R_{\log,f_0}).$$
b) The canonical principalization $\roman{cp}_1$ is an isomorphism in a neighborhood
of $p$, except under the subcase $a \geq c + 1$ in the case $i_s > i_o$.  In the
exceptional subcase, the strict transform of the exceptional divisor
$E_p$ for
$\roman{cp}_{1,1}$, obtained by blowing up $p$, does not appear in $R_{\log,f_1}$,
that is to say, the coefficient of the strict transform of the exceptional divisor
$E_p$ is 0, i.e.,
$$\nu_{E_p}(R_{\log,f_1}) = 0.$$
More generally, none of the irreducible components $E$ of the exceptional
divisor for
$\roman{cp}_1$ (i.e., none of the strict transforms of the exceptional
divisors for
$\roman{cp}_1, ... , \roman{cp}_l$), appear in
$R_{\log,f_1}$, i.e.,
$$\nu_E(R_{\log,f_1}) = 0.$$
\endproclaim

\demo{Proof}\enddemo a) Firstly we compute the coefficient
$\nu_{G_1}(R_{\log,f_0})$ of $G_1$ in the logarithmic ramification divisor
$R_{\log,f_0}$.
$$\align
\nu_{G_1}(R_{\log,f_0}) &= \nu_{G_1}\left(\frac{dy_1}{y_1} \wedge
dy_2/\frac{dx_1}{x_1} \wedge \frac{dx_2}{x_2}\right) \\
&= \nu_{G_1}\left(\frac{d(x_1^a)}{x_1^a} \wedge d\left(\Sigma_{i
\geq i_o > 0, j > 0} \alpha_{ij}x_1^ix_2^j\right)/\frac{dx_1}{x_1}
\wedge
\frac{dx_2}{x_2}\right) \\
&= \nu_{G_1}\left(\frac{d(x_1^a)}{x_1^a} \wedge d(x_1^{i_o}x_2)/\frac{dx_1}{x_1}
\wedge
\frac{dx_2}{x_2}\right) = i_o.\\
\endalign$$

Secondly we compute the coefficient of (the strict transform of) $G_1$ in
$R_{\log,f_1}$.

\vskip.1in

Subcase $a \geq \min\{i_o, i_s\}$: In this subcase, we have
$$\left\{\aligned
\nu_{G_1}(\frac{y_1}{y_2}) &= a - \min\{i_o, i_s\} \geq 0 \\
\nu_{G_1}(y_2) &= \min\{i_o, i_s\} > 0, \\
\endaligned\right.$$
which implies that the generic point of (the strict transform of) $G_1$ maps
under $f_1$ into the standard neighborhood of $q_1$ with a system of
regular parameters $(\frac{y_1}{y_2},y_2)$.  Therefore, we compute
$$\align
\nu_{G_1}(R_{\log,f_1}) &=
\nu_{G_1}\left(\frac{d(\frac{y_1}{y_2})}{(\frac{y_1}{y_2})} \wedge
\frac{dy_2}{y_2}/\frac{d(\frac{x_1}{x_2})}{(\frac{x_1}{x_2})} \wedge
\frac{dx_2}{x_2}\right) \\
&= \nu_{G_1}\left(\frac{dy_1}{y_1} \wedge
dy_2/\frac{dx_1}{x_1} \wedge \frac{dx_2}{x_2}\right) - \nu_{G_1}(y_2) \\
&= i_o - \min\{i_o,i_s\} < i_o.\\
\endalign$$

Subcase $a < \min\{i_o, i_s\}$: In this subcase, we have
$$\left\{\aligned
\nu_{G_1}(y_1) &= a > 0 \\
\nu_{G_1}(\frac{y_2}{y_1}) &= \min\{i_o, i_s\} - a > 0, \\
\endaligned\right.$$
which implies that the generic point of (the strict transform of) $G_1$ maps
under $f_1$ into the standard neighborhood of $q_1'$ with a system of
regular parameters $(y_1,\frac{y_2}{y_1})$.  Therefore, we compute
$$\align
\nu_{G_1}(R_{\log,f_1}) &=
\nu_{G_1}\left(\frac{dy_1}{y_1} \wedge
d(\frac{y_2}{y_1})/\frac{d(\frac{x_1}{x_2})}{(\frac{x_1}{x_2})} \wedge
\frac{dx_2}{x_2}\right) \\
&= \nu_{G_1}\left(\frac{dy_1}{y_1} \wedge
dy_2/\frac{dx_1}{x_1} \wedge \frac{dx_2}{x_2}\right) - \nu_{G_1}(y_1) \\
&= i_o - a < i_o.\\
\endalign$$
This completes the proof of the statement a).

\vskip.1in

b) We have only to consider the exceptional subcase $a \geq c + 1$ in the case $i_o >
i_1$.

Take the standard neighborhood of $q_1$ with
a system of regular parameters $(\frac{y_1}{y_2},y_2)$.  Then observing
$$\left\{\aligned
\nu_{G_1}(\frac{y_1}{y_2}) &= a - c > 0 \\
\nu_{G_1}(y_2) &= c > 0\\
\endaligned\right. \hskip.1in \& \hskip.1in \left\{\aligned
\nu_{E_p}(\frac{y_1}{y_2}) &= a - (c + 1) \geq 0 \\
\nu_{E_p}(y_2) &= c + 1 > 0,\\
\endaligned\right.$$
we conclude that the generic points of $G_1$ and $E_p$ maps into the standard
neighborhood of $q_1$.

Observe that the morphisms $f_0, \roman{bp}_1, \roman{cp}_{1,1}$ are all
toroidal with respect to the modified logarithmic structures obtained by
adding
$\{y_2 = 0\}$ and $\{x_2  = 0\}$ (and their strict transforms) to the original boundary
divisors.  The original logarithmic structures coincide with the modified ones
in neighborhoods of the generic points of $G_1$ and $E_p$ and in the standard
neighborhood of $q_1$.  Therefore, $f_{1,1}$ is toroidal in neighborhoods of the
generic points of $G_1$ and $E_p$ with respect to the original logarithmic structures,
and hence
$$\nu_{E_p}(R_{\log,f_1}) = 0$$
and
$$\nu_{G_1}(R_{\log,f_1}) = 0 < c = \nu_{G_1}(R_{\log,f_0}).$$
This proves statement a) and the first part of statement b).

The second part of statement b) follows inductively from the first part
of the statement b) and Claim 3.3.3.B.2 b).

\newpage

3.3.4) Subcase $2_p1_q0$: Neither $G_1$ nor $G_2$ maps onto $q$ (and hence
both
$G_1$ and $G_2$ would map onto $H_1$)

\vskip.1in

This subcase does not happen.  (Although this is again essentially due to
Abhyankar's lemma, we will not use it explicitly).  We will show by some
elementary calculation that, in this subcase,
$f|_{U_X}$ cannot be smooth and hence that we violate condition (c)
imposed on a morphism in the logarithmic category (cf. 1.2.1). 

\vskip.1in

First, we start with some ystems of regular parameters as described in 3.2. 
Then since $f^{-1}(D_Y) = D_X$ where $D_Y \cap U_q = H_1$ and $D_X
\cap U_p = G_1 \cup G_2$, we have
$$f^*y_1 = u \cdot x_1^ax_2^b \hskip.1in \text{for\ some\ unit}
\hskip.1in u \in \widehat{{\Cal O}_{X,p}}^{\times} \hskip.1in \text{with}
\hskip.1in a > 0, b > 0.$$ By replacing $x_1$ and $x_2$ with the ones
multiplied by some appropriate units, we may assume
$$f^*y_1 = x_1^ax_2^b \hskip.1in \text{with} \hskip.1in a > 0, b >
0.$$
On the other hand, since neither $G_1$ nor $G_2$ maps onto $q$, neither $x_1$
nor $x_2$ divides $f^*y_2$ (Note that $y_2$ can be chosen in any way,
irrelevant to the logarithmic structure, as long as $(y_1,y_2)$ form a system
of regular parameters.), in the Taylor expansion of $f^*y_2 = \Sigma \alpha_{ij}
x_1^ix_2^j$ we conclude that there exists $i > 0$ such that $\alpha_{i0} \neq
0$ and that there exists $j > 0$ such that $\alpha_{0j} \neq 0$.

We set
$$\align
i_o &:= \min \{i;\alpha_{i0} \neq 0\} > 0 \\
j_o &:= \min \{i;\alpha_{0j} \neq 0\} > 0. \\
\endalign$$

Finally we compute
$$\align
f^*(dy_1 \wedge dy_2) &= df^*y_1 \wedge df^*y_2 \\
&= (ax_1^{a-1}x_2^bdx_1 + x_1^abx_2^{b-1}dx_2)\\
&\wedge (\alpha_{i_o0}i_ox_1^{i_o-1}dx_1 + \Sigma_{i > i_o}
\alpha_{i0}ix_1^{i-1}dx_1 \\
& \hskip1.1in + d(\Sigma_{i > 0, j > 0} \alpha_{ij}x_1^ix_2^j) \\
&\hskip1.1in + \Sigma_{j > j_o}
\alpha_{0j}jx_2^{j-1}dx_2 + \alpha_{0j_o}j_ox_2^{j_o-1}dx_2)\\
&= (- b\alpha_{i_o0}i_ox_1^{i_o} - \Sigma_{i > i_o} b\alpha_{i0}ix_1^{i - i_o}
\cdot x_1^{i_o}\\
&\hskip1.1in + x_1x_2 \cdot h \\
&\hskip1.1in + a\alpha_{0j_o}j_ox_2^{j_o} + \Sigma_{j > j_o}
a\alpha_{0j}jx_2^{j - j_o}
\cdot x_2^{j_o}) \\
&\hskip1.1in \cdot x_1^{a-1}x_2^{b-1}dx_1
\wedge dx_2,\\
\endalign$$
noting that
$$(ax_1^{a-1}x_2^bdx_1 + x_1^abx_2^{b-1}dx_2) \wedge
d(\Sigma_{i > 0, j > 0} \alpha_{ij}x_1^ix_2^j)$$
is divisible by $x_1^ax_2^b \cdot dx_1 \wedge dx_2$ and can be
written as
$$x_1x_2 \cdot h \cdot x_1^{a-1}x_2^{b-1}dx_1
\wedge dx_2 \text{\ for\ some\ }h \in \widehat{{\Cal O}_{X,p}}.$$
But this implies that $f$ ramifies along
$$\{(- b\alpha_{i_o0}i_ox_1^{i_o} - \Sigma_{i > i_o} b\alpha_{i0}ix_1^{i - i_o}
\cdot x_1^{i_o} + x_1x_2 \cdot h + a\alpha_{0j_o}j_ox_2^{j_o} + \Sigma_{j >
j_o} a\alpha_{0j}jx_2^{j - j_o} \cdot x_2^{j_o}) = 0\}$$
other than possibly $\{x_1 = 0\}$ or $\{x_2 = 0\}$.  This violates condition
(c) $f|_{U_X = X - (\{x_1 = 0\} \cup \{x_2 = 0\})}$ being smooth imposed on
a morphism in the logarithmic category. 

\newpage

3.3.5) Subcase $2_p1_q1$: $G_1$ maps onto $q$ but $G_2$ does not map onto $q$

\vskip.1in

3.3.5.A) Coordinate Expression

\vskip.1in

First, we start with some systems of regular parameters as described in 3.2. 
Then since $f^{-1}(D_Y) = D_X$ where $D_Y \cap U_q = H_1$ and $D_X
\cap U_p = G_1 \cup G_2$, we have
$$f^*y_1 = u \cdot x_1^ax_2^b \hskip.1in \text{for\ some\ unit}
\hskip.1in u \in \widehat{{\Cal O}_{X,p}}^{\times} \hskip.1in \text{with}
\hskip.1in a > 0, b > 0.$$ By replacing $x_1$ and $x_2$ with the ones
multiplied by some appropriate units, we may assume
$$f^*y_1 = x_1^ax_2^b \hskip.1in \text{with} \hskip.1in a > 0, b >
0.$$
Set
$$f^*y_2 = \Sigma_{i > 0, j > 0, \det\left[\matrix
a & b \\
i & j \\
\endmatrix\right] = 0}\alpha_{ij}x_1^ix_2^j + \Sigma_{i > 0, j \geq 0,
\det\left[\matrix a & b \\
i & j \\
\endmatrix\right] \neq 0}\alpha_{ij}x_1^ix_2^j,$$
where no term of the form $x_2^j = x_1^0x_2^j\hskip.1in (i = 0, j > 0)$ appears,
since
$x_1$ has to divide $f^*y_2$ as $G_1 = \{x_1 = 0\}$ maps onto $q$, and where
there exists $i > 0$ with $\alpha_{i0} \neq 0$, since $x_2$ should not
divide
$f^*y_2$ as $G_2 = \{x_2 = 0\}$ does not map onto $q$.

Now since $f|_{U_X}$ must be smooth by condition (c) imposed on a morphism
in the logarithmic category (cf. 1.2.1) and by computing the
ramification \linebreak
$f^*(dy_1
\wedge dy_2)/dx_1 \wedge dx_2$ as in 3.3.2.A, we conclude that $f^*y_2$
must be of the form
$${\aligned
f^*y_2 &= \Sigma_{(i,j) \geq (i_s,j_s), \det\left[\matrix
a & b \\
i & j \\
\endmatrix\right] = 0}\alpha_{ij}x_1^ix_2^j \\
&+ \alpha_{i_o,0}x_1^{i_o} +
\Sigma_{i \geq i_o > 0, j \geq 0, (i,j) \neq (i_o,0),
\det\left[\matrix a & b \\
i & j \\
\endmatrix\right] \neq 0}\alpha_{ij}x_1^ix_2^j,\\
\endaligned} \hskip.1in \text{with} \hskip.1in \alpha_{i_o,0} \neq 0,$$
where
$$\align
i_0 &= \min\{i;\alpha_{i0} \neq 0\} \\
(i_s,j_s) &= \min\{(i,j);\alpha_{ij} \neq 0, \det\left[\matrix
a & b \\
i & j \\
\endmatrix\right] = 0\}.\\
\endalign$$

\vskip.1in

3.3.5.B) Descriptions of the canonical principalization $\roman{cp}_1$ and
the induced 

\hskip.3in morphism $f_1$

\vskip.1in

Case $i_o \leq i_s$ and $i_o \leq a$: In this case, we have
$$\left\{\aligned
f^*y_1 &= x_1^ax_2^b \\
f^*y_2 &= v \cdot x_1^{i_o} \\
\endaligned\right. \hskip.1in \text{with} \hskip.1in v \in
\widehat{{\Cal O}_{X,p}}^{\times} \hskip.1in
\text{unit},$$ which implies $f^{-1}(y_1,y_2) \cdot \widehat{{\Cal O}_{X,p}}
= (x_1^{i_o})$ is principal.  Therefore,
$\roman{cp}_1$ is an isomorphism in a neighborhood of $p$.

\vskip.1in

Case Otherwise, i.e., $i_o > i_s$ or $i_o > a$: We consider the 1st blowup
$\roman{cp}_{1,1}$ of the canonical principalization
$\roman{cp}_1$, which is a sequence of blowups $\roman{cp}_1 =
\roman{cp}_{1,1} \circ \cdot\cdot\cdot \circ \roman{cp}_{1,l}$ specified by the
canonical principalization algorithm.

\proclaim{Diagram 3.3.5.B.1}\endproclaim

\vskip2.5in

\proclaim{Claim 3.3.5.B.2} The canonical principalization $\roman{cp}_1$ is
an isomorphism in a neighborhood of $p$, except under the case $i_o > i_s$ or
$i_o > a$.  In the exceptional case, we claim:

a) The rational map $f_{1,1} = \roman{bp}^{-1} \circ f_0
\circ
\roman{cp}_{1,1}$ is well-defined (regular), that is to say, the ideal
$(f_0 \circ \roman{cp}_{1,1})^{-1}(y_1,y_2) \cdot {\Cal O}_{X_{1,1}}$ is
principal, except possibly at finitely many points on $E_p$.

b) At $p_1$, the morphism $f \circ \roman{cp}_{1,1}$ is in Subcase
$2_{p_1}1_q2$.  At $p_1'$, the morphism $f \circ \roman{cp}_{1,1}$ is in Subcase
$2_{p_1'}1_q1$.  At any other point $p_1'' (\neq p_1, p_1')$ on $E$, the
morphism
$f
\circ
\roman{cp}_{1,1}$ is in Subcase
$1_{p_1''}1_q1$.
\endproclaim
 \demo{Proof}\enddemo Statements a) and b) are obvious.

\vskip.1in

3.3.5.C) Conclusion on the behavior of the logarithmic ramification divisor

\vskip.1in

\proclaim{Claim 3.3.5.C.1} a) The coefficient of $G_1$ in the logarithmic
ramification divisor strictly drops, while the coefficient of $G_2$ in the
logarithmic ramification divisor remains 0.  More precisely,
$$\nu_{G_1}(R_{\log,f_1}) = \left\{\aligned
& i_o - \min\{i_o,i_s\} \hskip.1in \text{if} \hskip.1in a \geq \min\{i_o,i_s\}
\\ 
& i_o - a \hskip.67in \text{if} \hskip.1in a < \min\{i_o,i_s\}
\\ 
\endaligned\right\} < i_o = \nu_{G_1}(R_{\log,f_0})$$
and 
$$\nu_{G_2}(R_{\log,f_1}) = \nu_{G_2}(R_{\log,f_0}) = 0.$$
b) The canonical principalization is an isomorphism in a
neighborhood of $p$, except under the case $i_o > i_s$ or $i_o > a$.  In the
exceptional case, the coefficient in $R_{\log,f_1}$ of the strict transform of the
eceptional divisor
$E_p$ for $\roman{cp}_{1,1}$, obtained by blowing up $p$, is strictly
smaller than the maximum of the coefficients of irreducible components in
$R_{\log,f_0}$ in a neighborhood of $p$.  More precisely, 
$$\nu_{E_p}(R_{\log,f_1}) = \left\{\aligned
& i_o - \min\{i_o,i_s + j_s\} \hskip.04in \text{if} \hskip.1in a + b \geq
\min\{i_o,i_s + j_s\}
\\ 
& i_o - (a + b) \hskip.56in \text{if} \hskip.1in a + b < \min\{i_o,i_s +
j_s\}
\\ 
\endaligned\right\} < i_o = \nu_{G_1}(R_{\log,f_0}).$$
More generally, the coefficient in $R_{\log,f_1}$, of any irreducible
component $E$ of the exceptional divisor for
$\roman{cp}_1$, is strictly smaller than the maximum of the coefficients of irreducible components in
$R_{\log,f_0}$ in a neighborhood of $p$, i.e.,
$$\nu_E(R_{\log,f_1}) < i_o = \nu_{G_1}(R_{\log,f_0}).$$
\endproclaim

\demo{Proof}\enddemo a) Firstly we compute the coefficient
$\nu_{G_1}(R_{\log,f_0})$ of $G_1$ in the logarithmic ramification divisor
$R_{\log,f_0}$:
$$\align
\nu_{G_1}(R_{\log,f_0}) &= \nu_{G_1}\left(\frac{dy_1}{y_1} \wedge
dy_2/\frac{dx_1}{x_1} \wedge \frac{dx_2}{x_2}\right) \\
&= \nu_{G_1}\left(\frac{d(x_1^ax_2^b)}{x_1^ax_2^b} \wedge d\left(\Sigma_{i
\geq i_o > 0, j \geq 0, \det\left[
\matrix
a & b \\
i & j \\
\endmatrix\right] \neq 0} \alpha_{ij}x_1^ix_2^j\right)/\frac{dx_1}{x_1}
\wedge
\frac{dx_2}{x_2}\right) \\
&= \nu_{G_1}\left(\frac{d(x_1^ax_2^b)}{x_1^ax_2^b} \wedge d(x_1^{i_o})/\frac{dx_1}{x_1}
\wedge
\frac{dx_2}{x_2}\right) = i_o.\\
\endalign$$

Secondly we compute the coefficient of (the strict transform of) $G_1$ in
$R_{\log,f_1}$.

\vskip.1in

Subcase $a \geq \min\{i_o, i_s\}$: In this subcase, we have
$$\left\{\aligned
\nu_{G_1}(\frac{y_1}{y_2}) &= a - \min\{i_o, i_s\} \geq 0 \\
\nu_{G_1}(y_2) &= \min\{i_o, i_s\} > 0, \\
\endaligned\right.$$
which implies that the generic point of (the strict transform of) $G_1$ maps
under $f_1$ into the standard neighborhood of $q_1$ with a system of
regular parameters $(\frac{y_1}{y_2},y_2)$.  Therefore, we compute
$$\align
\nu_{G_1}(R_{\log,f_1}) &=
\nu_{G_1}\left(\frac{d(\frac{y_1}{y_2})}{(\frac{y_1}{y_2})} \wedge
\frac{dy_2}{y_2}/\frac{d(\frac{x_1}{x_2})}{(\frac{x_1}{x_2})} \wedge
\frac{dx_2}{x_2}\right) \\
&= \nu_{G_1}\left(\frac{dy_1}{y_1} \wedge
dy_2/\frac{dx_1}{x_1} \wedge \frac{dx_2}{x_2}\right) - \nu_{G_1}(y_2) \\
&= i_o - \min\{i_o,i_s\} < i_o.\\
\endalign$$

Subcase $a < \min\{i_o, i_s\}$: In this subcase, we have
$$\left\{\aligned
\nu_{G_1}(y_1) &= a > 0 \\
\nu_{G_1}(\frac{y_2}{y_1}) &= \min\{i_o, i_s\} - a > 0, \\
\endaligned\right.$$
which implies that the generic point of (the strict transform of) $G_1$ maps
under $f_1$ into the standard neighborhood of $q_1'$ with a system of
regular parameters $(y_1,\frac{y_2}{y_1})$.  Therefore, we compute
$$\align
\nu_{G_1}(R_{\log,f_1}) &=
\nu_{G_1}\left(\frac{dy_1}{y_1} \wedge
d(\frac{y_2}{y_1})/\frac{d(\frac{x_1}{x_2})}{(\frac{x_1}{x_2})} \wedge
\frac{dx_2}{x_2}\right) \\
&= \nu_{G_1}\left(\frac{dy_1}{y_1} \wedge
dy_2/\frac{dx_1}{x_1} \wedge \frac{dx_2}{x_2}\right) - \nu_{G_1}(y_1) \\
&= i_o - a < i_o.\\
\endalign$$
This completes the proof of statement a).

\vskip.1in

b) We have only to consider the exceptional case $i_o > i_s$ or $i_o > a$.

We compute the coefficient of (the strict transform of) the exceptional
divisor $E_p$ for $\roman{cp}_{1,1}$, obtained by blowing up $p$, in the
logarithmic ramification divisor $R_{\log,f_1}$.

\vskip.1in

Subcase $a + b \geq \min\{i_o, i_s + j_s\}$: In this subcase, we have
$$\left\{\aligned
\nu_{E_p}(\frac{y_1}{y_2}) &= a + b - \min\{i_o, i_s + j_s\} \geq 0 \\
\nu_{E_p}(y_2) &= \min\{i_o, i_s + j_s\} > 0, \\
\endaligned\right.$$
which implies that the generic point of (the strict transform of) $E_p$ maps
under $f_1$ into the standard neighborhood of $q_1$ with a system of
regular parameters $(\frac{y_1}{y_2},y_2)$.  Therefore, we compute
$$\align
\nu_{E_p}(R_{\log,f_1}) &=
\nu_{E_p}\left(\frac{d(\frac{y_1}{y_2})}{(\frac{y_1}{y_2})} \wedge
\frac{dy_2}{y_2}/\frac{d(\frac{x_1}{x_2})}{(\frac{x_1}{x_2})} \wedge
\frac{dx_2}{x_2}\right) \\
&= \nu_{E_p}\left(\frac{dy_1}{y_1} \wedge
dy_2/\frac{dx_1}{x_1} \wedge \frac{dx_2}{x_2}\right) - \nu_{E_p}(y_2) \\
&= i_o - \min\{i_o,i_s + j_s\} < i_o.\\
\endalign$$

\vskip.1in

Subcase $a + b < \min\{i_o, i_s + j_s\}$: In this subcase, we have
$$\left\{\aligned
\nu_{E_p}(y_1) &= a + b > 0 \\
\nu_{E_p}(\frac{y_2}{y_1}) &= \min\{i_o, i_s + j_s\} - (a + b) > 0, \\
\endaligned\right.$$
which implies that the generic point of (the strict transform of) $E_p$ maps
under $f_1$ into the standard neighborhood of $q_1'$ with a system of
regular parameters $(y_1,\frac{y_2}{y_1})$.  Therefore, we compute
$$\align
\nu_{E_p}(R_{\log,f_1}) &=
\nu_{E_p}\left(\frac{dy_1}{y_1} \wedge
d(\frac{y_2}{y_1})/\frac{d(\frac{x_1}{x_2})}{(\frac{x_1}{x_2})} \wedge
\frac{dx_2}{x_2}\right) \\
&= \nu_{E_p}\left(\frac{dy_1}{y_1} \wedge
dy_2/\frac{dx_1}{x_1} \wedge \frac{dx_2}{x_2}\right) - \nu_{E_p}(y_1) \\
&= i_o - (a + b) < i_o.\\
\endalign$$

This completes the proof of the first part of statement b).

\vskip.1in

In order to verify the second part of statement b), we look
at the points where $(f \circ \roman{cp}_{1,1})^{-1}(y_1,y_2) \cdot {\Cal
O}_{X_{1,1}}$ is (possibly) not principal.

At $p_1$, the morphism $f \circ \roman{cp}_{1,1}$ is in Subcase
$2_{p_1}1_q2$ with
$$\align
\nu_{G_1}(R_{\log,f \circ \roman{cp}_{1,1}}) &= i_o \\
\nu_{E_p}(R_{\log,f \circ \roman{cp}_{1,1}}) &= i_o. \\
\endalign$$
Therefore, by Claim 3.3.6.C.1 b) we have for any irreducible component $E$,
other than (the strict transform of) $E_p$, which is exceptional for
$\roman{cp}_1$ and which maps onto $p_1$ under $\roman{cp}_{1,2} \circ
\cdot\cdot\cdot
\circ
\roman{cp}_{1,l}$
$$\nu_E(R_{\log,f_1}) < \max\{\nu_{G_1}(R_{\log,f \circ
\roman{cp}_{1,1}}),\nu_{E_p}(R_{\log,f \circ \roman{cp}_{1,1}})\} = i_o =
\nu_{G_1}(R_{\log,f_0}).$$

At $p_1'$, the morphism $f \circ \roman{cp}_{1,1}$ is in Subcase
$2_{p_1'}1_q1$ with
$$\nu_{E_p}(R_{\log,f \circ \roman{cp}_{1,1}}) = i_o.$$
Therefore, by the first part of statement b) and by induction on the
length of the sequence of blowups of points for $\roman{cp}_1$, we have for
any irreducible component $E$, other than (the strict transform of) $E_p$,
which is exceptional for
$\roman{cp}_1$ and which maps onto $p_1'$ under $\roman{cp}_{1,2} \circ
\cdot\cdot\cdot
\circ
\roman{cp}_{1,l}$
$$\nu_E(R_{\log,f_1}) < \nu_{E_p}(R_{\log,f \circ \roman{cp}_{1,1}}) = i_o =
\nu_{G_1}(R_{\log,f_0}).$$

At $p_1'' \in E_p$, other than $p_1$ or $p_1'$, the morphism $(f \circ
\roman{cp}_{1,1})$ is in Subcase $1_{p_1''}1_q1$.  Therefore, by Claim
3.3.3.C.1 b) we have for any irreducible component $E$,
other than (the strict transform of) $E_p$, which is exceptional for
$\roman{cp}_1$ and which maps onto $p_1''$ under $\roman{cp}_{1,2} \circ
\cdot\cdot\cdot
\circ
\roman{cp}_{1,l}$
$$\nu_E(R_{\log,f_1}) = 0 < \nu_{E_p}(R_{\log,f \circ \roman{cp}_{1,1}}) =
i_o =
\nu_{G_1}(R_{\log,f_0}).$$

\newpage

3.3.6) Subcase $2_p1_q2$: Both $G_1$ and $G_2$ map onto $q$

\vskip.1in

3.3.6.A) Coordinate Expression

\vskip.1in

First, we start with some systems of regular parameters as described in
3.2.  Then since $f^{-1}(D_Y) = D_X$ where $D_Y \cap U_q = H_1$ and $D_X
\cap U_p = G_1 \cup G_2$, we have
$$f^*y_1 = u \cdot x_1^ax_2^b \hskip.1in \text{for\ some\ unit}
\hskip.1in u \in \widehat{{\Cal O}_{X,p}}^{\times} \hskip.1in \text{with}
\hskip.1in a > 0, b > 0.$$ By replacing $x_1$ and $x_2$ with the ones
multiplied by some appropriate units, we may assume
$$f^*y_1 = x_1^ax_2^b \hskip.1in \text{with} \hskip.1in a > 0, b >
0.$$
Set
$$f^*y_2 = \Sigma_{i > 0, j > 0, \det\left[\matrix
a & b \\
i & j \\
\endmatrix\right] = 0}\alpha_{ij}x_1^ix_2^j + \Sigma_{i > 0, j > 0,
\det\left[\matrix a & b \\
i & j \\
\endmatrix\right] \neq 0}\alpha_{ij}x_1^ix_2^j,$$
where $x_1x_2$ divides $f^*y_2$, since both $G_1$ and $G_2$ maps onto $q$.  

Now since $f|_{U_X}$ must be smooth by condition (c) imposed on a morphism
in the logarithmic category (cf. 1.2.1) and by computing the
ramification \linebreak
$f^*(dy_1
\wedge dy_2)/dx_1 \wedge dx_2$ as in 3.3.2.A, we conclude that $f^*y_2$
must be of the form
$${\aligned
f^*y_2 &= \Sigma_{(i,j) \geq (i_s,j_s), \det\left[\matrix
a & b \\
i & j \\
\endmatrix\right] = 0}\alpha_{ij}x_1^ix_2^j\\
&+ \alpha_{i_o,j_o}x_1^{i_o}x_2^{j_o} + \Sigma_{i \geq i_o > 0, j \geq j_o >
0, (i,j) \neq (i_o,j_o), \det\left[\matrix a & b \\
i & j \\
\endmatrix\right] \neq 0}\alpha_{ij}x_1^ix_2^j,\\
\endaligned} \hskip.1in \text{with} \hskip.1in \alpha_{i_o,j_o} \neq 0,$$
where
$$\align
i_o &= \min\{i;\alpha_{ij} \neq 0, \det\left[\matrix
a & b \\
i & j \\
\endmatrix\right] \neq 0\} \\
j_o &= \min\{j;\alpha_{ij} \neq 0, \det\left[\matrix
a & b \\
i & j \\
\endmatrix\right] \neq 0\} \\
&\text{\ with\ }\det\left[\matrix
a & b \\
i_o & j_o \\
\endmatrix\right] \neq 0, \text{\ and\ } \\
(i_s,j_s) &= \min\{(i,j);\alpha_{ij} \neq 0, \det\left[\matrix
a & b \\
i & j \\
\endmatrix\right] = 0\}.\\
\endalign$$

\vskip.1in

3.3.6.B) Descriptions of the canonical principalization $\roman{cp}_1$ and
the induced 

\hskip.3in morphism $f_1$

\vskip.1in

Case $\left\{\aligned
&i_o \geq i_s \\
&j_o \geq j_s \\
\endaligned\right\}$
or
$\left\{\aligned
&i_o \leq i_s \\
&j_o \leq j_s \\
\endaligned \hskip.1in \& \hskip.1in
\aligned 
&i_o \leq a \\
&j_o \leq b \\
\endaligned
\right\}$
or
$\left\{\aligned
&i_o \leq i_s \\
&j_o \leq j_s \\
\endaligned \hskip.1in \& \hskip.1in
\aligned 
&i_o \geq a \\
&j_o \geq b \\
\endaligned
\right\}$: 

\vskip.1in

In this case, the ideal 
$$\align
f^{-1}(y_1,y_2) \cdot \widehat{{\Cal O}_{X,p}} &= (x_1^{i_s}x_2^{j_s})
\hskip.1in
\text{or} \hskip.1in (x_1^ax_2^b), \\
&\text{or}\\
&(x_1^{i_o}x_2^{j_o}), \\
&\text{or}\\
&(x_1^ax_2^b), \hskip.1in \text{respectively}, \\
\endalign$$
is principal.  Therefore,
$\roman{cp}_1$ is an isomorphism in a neighborhood of $p$.

\vskip.1in

Case Otherwise:

\vskip.1in

We consider the 1st blowup $\roman{cp}_{1,1}$ of the canonical principalization
$\roman{cp}_1$, which is a sequence of blowups $\roman{cp}_1 =
\roman{cp}_{1,1} \circ \cdot\cdot\cdot \circ \roman{cp}_{1,l}$ specified by the
canonical principalization algorithm.

\proclaim{Diagram 3.3.6.B.1}\endproclaim

\vskip2.5in

\proclaim{Claim 3.3.6.B.2} The canonical principalization is an isomorphism
in a neighborhood of $p$, except for the case where $f^{-1}(y_1,y_2) \cdot
\widehat{{\Cal O}_{X,p}}$ is not principal.  In the exceptional case, we
claim:

a) The rational map $f_{1,1} = \roman{bp}^{-1} \circ f_0
\circ
\roman{cp}_{1,1}$ is well-defined (regular), that is to say, the ideal
$(f_0 \circ \roman{cp}_{1,1})^{-1}(y_1,y_2) \cdot {\Cal O}_{X_{1,1}}$ is
principal, except possibly at finitely many points on $E_p$.

b) At $p_1$ (resp. $p_1'$), the morphism
$f
\circ
\roman{cp}_{1,1}$ is in Subcase
$2_{p_1}1_q1$ (resp. $2_{p_1'}1_q1$).  At any other point $p_1''$ on $E_p$, the
morphism
$f
\circ
\roman{cp}_{1,1}$ is in Subcase
$1_{p_1''}1_q1$.
\endproclaim
\demo{Proof}\enddemo Statements a) and b) are obvious.

\vskip.1in

3.3.6.C) Conclusion on the behavior of the logarithmic ramification divisor

\vskip.1in

\proclaim{Claim 3.3.6.C.1} a) The coefficients of $G_1$ and $G_2$ in the
logarithmic ramification divisor strictly drop.  More precisely,
$$\nu_{G_1}(R_{\log,f_1}) = \left\{\aligned
& i_o - \min\{i_o,i_s\} \hskip.1in \text{if} \hskip.1in a \geq \min\{i_o,i_s\}
\\ 
& i_o - a \hskip.67in \text{if} \hskip.1in a < \min\{i_o,i_s\}
\\ 
\endaligned\right\} < i_o = \nu_{G_1}(R_{\log,f_0})$$
$$\nu_{G_2}(R_{\log,f_1}) = \left\{\aligned
& j_o - \min\{j_o,j_s\} \hskip.1in \text{if} \hskip.1in b \geq
\min\{j_o,j_s\}
\\ 
& j_o - b \hskip.7in \text{if} \hskip.1in b < \min\{j_o,j_s\}
\\ 
\endaligned\right\} < j_o = \nu_{G_2}(R_{\log,f_0}).$$

b) The canonical principalization is an isomorphism in a
neighborhood of $p$, except for the case when $f^{-1}(y_1,y_2) \cdot
\widehat{{\Cal O}_{X,p}}$ is not principal.  In the exceptional case where
$f^{-1}(y_1,y_2) \cdot \widehat{{\Cal O}_{X,p}}$ is not principal, the
coefficient in
$R_{\log,f_1}$ of the strict transform of the eceptional divisor
$E_p$ for $\roman{cp}_{1,1}$, obtained by blowing up $p$, is strictly
smaller than the maximum of the coefficients of the irreducible components in
$R_{\log,f_0}$ in a neighborhood of $p$.  More precisely, 
$$\align
\nu_{E_p}(R_{\log,f_1}) &= \left\{\aligned
& i_o + j_o - \min\{i_o + j_o,i_1 + j_1\} \hskip.1in \text{if} \hskip.1in a
+ b
\geq
\min\{i_o + j_o,i_1 + j_1\}
\\ 
& i_o + j_o - (a + b) \hskip.92in \text{if} \hskip.1in a + b < \min\{i_o +
j_o,i_1 + j_1\}
\\ 
\endaligned\right\} \\
&< \max\{i_o,j_o\} = \max\{\nu_{G_1}(R_{\log,f_0}), \nu_{G_2}(R_{\log,f_0})\}.\\
\endalign$$
More generally, the coefficient in $R_{\log,f_1}$, of any irreducible
component $E$ of the exceptional divisor for
$\roman{cp}_1$, is strictly smaller than the maximum of the coefficients of
the irreducible components in
$R_{\log,f_0}$ in a neighborhood of $p$, i.e.,
$$\nu_E(R_{\log,f_1}) < \max\{\nu_{G_1}(R_{\log,f_0}),
\nu_{G_1}(R_{\log,f_0})\} = \max\{i_o,j_o\}.$$
\endproclaim

\demo{Proof}\enddemo a) Firstly we compute the coefficient
$\nu_{G_1}(R_{\log,f_0})$ of $G_1$ in the logarithmic ramification divisor
$R_{\log,f_0}$.
$$\align
\nu_{G_1}(R_{\log,f_0}) &= \nu_{G_1}\left(\frac{dy_1}{y_1} \wedge
dy_2/\frac{dx_1}{x_1} \wedge \frac{dx_2}{x_2}\right) \\
&= \nu_{G_1}\left(\frac{d(x_1^ax_2^b)}{x_1^ax_2^b} \wedge d\left(\Sigma_{i
\geq i_o > 0, j \geq j_o > 0, \det\left[
\matrix
a & b \\
i & j \\
\endmatrix\right] \neq 0} \alpha_{ij}x_1^ix_2^j\right)/\frac{dx_1}{x_1}
\wedge
\frac{dx_2}{x_2}\right) \\
&= \nu_{G_1}\left(\frac{d(x_1^ax_2^b)}{x_1^ax_2^b} \wedge
d(x_1^{i_o}x_2^{j_o})/\frac{dx_1}{x_1}
\wedge
\frac{dx_2}{x_2}\right) = i_o.\\
\endalign$$
Symmetrically, we compute the coefficient
$\nu_{G_2}(R_{\log,f_0})$ of $G_2$ in the logarithmic ramification divisor
$R_{\log,f_0}$ to be
$$\nu_{G_2}(R_{\log,f_0}) = j_o.$$

Secondly we compute the coefficient of (the strict transform of) $G_1$ in
$R_{\log,f_1}$.

\vskip.1in

Subcase $a \geq \min\{i_o, i_s\}$: In this subcase, we have
$$\left\{\aligned
\nu_{G_1}(\frac{y_1}{y_2}) &= a - \min\{i_o, i_s\} \geq 0 \\
\nu_{G_1}(y_2) &= \min\{i_o, i_s\} > 0, \\
\endaligned\right.$$
which implies that the generic point of (the strict transform of) $G_1$ maps
under $f_1$ into the standard neighborhood of $q_1$ with a system of
regular parameters $(\frac{y_1}{y_2},y_2)$.  Therefore, we compute
$$\align
\nu_{G_1}(R_{\log,f_1}) &=
\nu_{G_1}\left(\frac{d(\frac{y_1}{y_2})}{(\frac{y_1}{y_2})} \wedge
\frac{dy_2}{y_2}/\frac{d(\frac{x_1}{x_2})}{(\frac{x_1}{x_2})} \wedge
\frac{dx_2}{x_2}\right) \\
&= \nu_{G_1}\left(\frac{dy_1}{y_1} \wedge
dy_2/\frac{dx_1}{x_1} \wedge \frac{dx_2}{x_2}\right) - \nu_{G_1}(y_2) \\
&= i_o - \min\{i_o,i_s\} < i_o.\\
\endalign$$

Subcase $a < \min\{i_o, i_s\}$: In this subcase, we have
$$\left\{\aligned
\nu_{G_1}(y_1) &= a > 0 \\
\nu_{G_1}(\frac{y_2}{y_1}) &= \min\{i_o, i_s\} - a > 0, \\
\endaligned\right.$$
which implies that the generic point of (the strict transform of) $G_1$ maps
under $f_1$ into the standard neighborhood of $q_1'$ with a system of
regular parameters $(y_1,\frac{y_2}{y_1})$.  Therefore, we compute
$$\align
\nu_{G_1}(R_{\log,f_1}) &=
\nu_{G_1}\left(\frac{dy_1}{y_1} \wedge
d(\frac{y_2}{y_1})/\frac{d(\frac{x_1}{x_2})}{(\frac{x_1}{x_2})} \wedge
\frac{dx_2}{x_2}\right) \\
&= \nu_{G_1}\left(\frac{dy_1}{y_1} \wedge
dy_2/\frac{dx_1}{x_1} \wedge \frac{dx_2}{x_2}\right) - \nu_{G_1}(y_1) \\
&= i_o - a < i_o.\\
\endalign$$
Symmetrically, we compute
$$\nu_{G_2}(R_{\log,f_1}) = \left\{\aligned
& j_o - \min\{j_o,j_s\} \hskip.1in \text{if} \hskip.1in b \geq
\min\{j_o,j_s\}
\\ 
& j_o - b \hskip.71in \text{if} \hskip.1in b < \min\{j_o,j_s\}
\\ 
\endaligned\right\} < j_o = \nu_{G_2}(R_{\log,f_0}).$$
This proves statement a).

\vskip.1in

b) We compute the coefficient of (the strict transform of) the exceptional
divisor $E_p$ for $\roman{cp}_{1,1}$, obtained by blowing up $p$, in the
logarithmic ramification divisor $R_{\log,f_1}$.

\vskip.1in

Subcase $a + b \geq \min\{i_o + j_o, i_s + j_s\}$: In this subcase, we have
$$\left\{\aligned
\nu_{E_p}(\frac{y_1}{y_2}) &= a + b - \min\{i_o + j_o, i_s + j_s\} \geq 0 \\
\nu_{E_p}(y_2) &= \min\{i_o + j_o, i_s + j_s\} > 0, \\
\endaligned\right.$$
which implies that the generic point of (the strict transform of) $E_p$ maps
under $f_1$ into the standard neighborhood of $q_1$ with a system of
regular parameters $(\frac{y_1}{y_2},y_2)$.  Therefore, we compute
$$\align
\nu_{E_p}(R_{\log,f_1}) &=
\nu_{E_p}\left(\frac{d(\frac{y_1}{y_2})}{(\frac{y_1}{y_2})} \wedge
\frac{dy_2}{y_2}/\frac{d(\frac{x_1}{x_2})}{(\frac{x_1}{x_2})} \wedge
\frac{dx_2}{x_2}\right) \\
&= \nu_{E_p}\left(\frac{dy_1}{y_1} \wedge
dy_2/\frac{dx_1}{x_1} \wedge \frac{dx_2}{x_2}\right) - \nu_{E_p}(y_2) \\
&= i_o + j_o - \min\{i_o + j_o,i_s + j_s\}.\\
\endalign$$
Now we claim that
$$i_o + j_o - \min\{i_o + j_o,i_s + j_s\} < \max\{i_o,j_o\}$$
unless $f^{-1}(y_1,y_2) \cdot \widehat{{\Cal O}_{X,p}}$ is already principal.

In fact, if
$$i_o + j_o \leq i_s + j_s,$$
then
$$i_o + j_o - \min\{i_o + j_o, i_s + j_s\} = 0 < \max\{i_o,j_o\}.$$
Therefore, we may assume
$$i_o + j_o > i_s + j_s,$$
in which case we have
$$\align
i_o + j_o &- \min\{i_o + j_o, i_s + j_s\} = (i_o + j_o) - (i_s + j_s) \\
&= (i_o - i_s) + (j_o - j_s) \\
&=\left\{\aligned
&(i_o - i_s) - (j_s - j_o) < i_o \hskip.2in \text{\ if\ }j_s \geq j_o \\
&(j_o - j_s) - (i_s - i_o) < j_o \hskip.2in \text{\ if\ }i_s \geq i_o \\
& \text{or} \\
&f^{-1}(y_1,y_2) \cdot \widehat{{\Cal O}_{X,p}} = (x_1^ax_2^b,
x_1^{i_s}x_2^{j_s}) = (x_1^{\min\{a,i_s\}}x_2^{\min\{b,j_s\}})\\
&\hskip1.06in \text{\ is\ principal\ }
\text{\ if\ }i_s < i_o
\hskip.05in \& \hskip.05in j_s < j_o.\\
\endaligned\right.\\
\endalign$$

\vskip.1in

Subcase $a + b < \min\{i_o + j_o, i_s + j_s\}$: In this subcase, we have
$$\left\{\aligned
\nu_{E_p}(y_1) &= a + b > 0 \\
\nu_{E_p}(\frac{y_2}{y_1}) &= \min\{i_o + j_o, i_s + j_s\} - (a + b) > 0, \\
\endaligned\right.$$
which implies that the generic point of (the strict transform of) $E_p$ maps
under $f_1$ into the standard neighborhood of $q_1'$ with a system of
regular parameters $(y_1,\frac{y_2}{y_1})$.  Therefore, we compute
$$\align
\nu_{E_p}(R_{\log,f_1}) &=
\nu_{E_p}\left(\frac{dy_1}{y_1} \wedge
d(\frac{y_2}{y_1})/\frac{d(\frac{x_1}{x_2})}{(\frac{x_1}{x_2})} \wedge
\frac{dx_2}{x_2}\right) \\
&= \nu_{E_p}\left(\frac{dy_1}{y_1} \wedge
dy_2/\frac{dx_1}{x_1} \wedge \frac{dx_2}{x_2}\right) - \nu_{E_p}(y_1) \\
&= i_o + j_o - (a + b).\\
\endalign$$

Now we claim that
$$i_o + j_o - (a + b) < \max\{i_o,j_o\}$$
unless $f^{-1}(y_1,y_2) \cdot \widehat{{\Cal O}_{X,p}}$ is already principal.

In fact, we have
$$\align
i_o + j_o - (a + b) &= (i_o - a) + (j_o - b) \\
&=\left\{\aligned
&(i_o - a) - (b - j_o) < i_o \hskip.2in \text{\ if\ }b \geq j_o \\
&(j_o - b) - (a - i_o) < j_o \hskip.2in \text{\ if\ }a \geq i_o \\
&\text{or}\\
&f^{-1}(y_1,y_2) \cdot \widehat{{\Cal O}_{X,p}} = (x_1^ax_2^b) \\
&\hskip.95in \text{\ principal\ } \text{\ if\
}a < i_o
\hskip.05in \& \hskip.05in b < j_o\\
&(\text{Note\ that\ } a < i_s \hskip.05in \& \hskip.05in b < j_s \\
&\text{\
by\ the\ subcase\ assumption\ } a + b < \min\{i_o + j_o, i_s + j_s\}.)\\
\endaligned\right.\\
\endalign$$
This proves the first part of statement b).

\vskip.1in

It remains to prove the second part of statement b).

\vskip.1in

Suppose $p_k \in D_{X_{1,k}} \subset X_{1,k}$ is a closed point sitting over
\linebreak
$p = p_0
\in D_X = D_{X_{1,0}} \subset X = X_{1,0}$, i.e.,
$$\roman{cp}_{1,1} \circ \cdot\cdot\cdot \circ \roman{cp}_{1,k}(p_k) = p.$$
We set
$$g_k = f \circ (\roman{cp}_{1,1} \circ \cdot\cdot\cdot \circ \roman{cp}_{1,k}).$$
Suppose that $g_k^{-1}(y_1,y_2) \cdot \widehat{{\Cal O}_{X_{1,k},p_k}}$ is
not principal, and hence that we blowup $p_k$ to obtain the exceptional
divisor
$E_{p_k}$ for
$\roman{cp}_{1,k+1}$, the $(k+1)$-th stage of the canonical principalization
$\roman{cp}_1 = \roman{cp}_{1,1} \circ \cdot\cdot\cdot \circ \roman{cp}_{1,k} \circ
\roman{cp}_{1,k+1} \circ \cdot\cdot\cdot \circ \roman{cp}_{1,l}$.

\vskip.1in

At $p_k$, the morphism $g_k$ is either in Subcase $1_{p_k}1_q1$ or in Subcase
$2_{p_k}1_q2$.

\vskip.1in

If $g_k$ is in Subcase $1_{p_k}1_q1$ at $p_k$, then by Claim 3.3.3.C.1 we
have
$$\nu_E(R_{\log,f_1}) = 0 < \max\{i_o,j_o\}$$
for (the strict transform of) any exceptional divisor $E$ for $\roman{cp}_{1,k+1}
\circ \cdot\cdot\cdot \circ \roman{cp}_{1,l}$ which maps onto $p_k$.

\vskip.1in

In order to analyze the case where $g_k$ is in Subcase $2_{p_k}1_q2$, we prove the
following inductive lemma.

\proclaim{Lemma 3.3.6.C.2} Let $p_k \in D_{X_{1,k}} \subset X_{1,k}$ be a
point, sitting over $p$, where
$g_k$ is in Subcase
$2_{p_k}1_q2$.  Assume that the ideal $g_k^{-1}(y_1,y_2) \cdot \widehat{{\Cal
O}_{X_{1,k},p_k}}$ is not principal.  

Suppose that we have the following two conditions: 

$(i)_k$ There exists a system of regular parameters
$(x_{1,k},x_{2,k})$ of $\widehat{{\Cal O}_{X_{1,k},p_k}}$, compatible with the
logarithmic structure, such that
$$(*)_k\left\{\aligned
g_k^*y_1 &= x_{1,k}^{a_k}x_{2,k}^{b_k} \hskip.1in \text{with} \hskip.1in a_k > 0, b_k >
0 \\
g_k^*y_2 &= \Sigma_{(i,j) \geq (i_{s,k},j_{s,k}), \det\left[\matrix
a_k & b_k \\
i & j \\
\endmatrix\right] = 0}\alpha_{ij}x_{1,k}^ix_{2,k}^j\\
&+ \alpha_{i_{o,k},j_{o,k}}x_{1,k}^{i_{o,k}}x_{2,k}^{j_{o,k}} + \Sigma_{i \geq i_{o,k}
> 0, j
\geq j_{o,k} > 0, (i,j) \neq (i_{o,k},j_{o,k}), \det\left[\matrix a_k & b_k
\\ i & j \\
\endmatrix\right] \neq 0}\alpha_{ij}x_{1,k}^ix_{2,k}^j\\
&\hskip.2in \text{with} \hskip.1in \alpha_{i_{o,k}, j_{o,k}} \neq 0, \\
\endaligned\right.$$
where
$$\align
i_{o,k} &= \min\{i;\alpha_{ij} \neq 0, \det\left[\matrix
a_k & b_k \\
i & j \\
\endmatrix\right] \neq 0\} \\
j_{o,k} &= \min\{j;\alpha_{ij} \neq 0, \det\left[\matrix
a_k & b_k \\
i & j \\
\endmatrix\right] \neq 0\} \\
&\text{\ with\ }\det\left[\matrix
a_k & b_k \\
i_{o,k} & j_{o,k} \\
\endmatrix\right] \neq 0, \text{\ and\ } \\
(i_{s,k},j_{s,k}) &= \min\{(i,j);\alpha_{ij} \neq 0, \det\left[\matrix
a_k & b_k \\
i & j \\
\endmatrix\right] = 0\}.\\
\endalign$$

\newpage

$(ii)_k$ One of the following four holds:
$$(\alpha)_k \left\{\aligned
i_{o,k} &> i_{s,k} \\
j_{o,k} &< j_{s,k} \\
\endaligned\right\} \hskip.1in \text{with} \hskip.1in i_{o,k} - i_{s,k} < i_o
\hskip.1in \text{and} \hskip.1in i_{o,k} - a_k < i_o$$

$$(\beta)_k \left\{\aligned
i_{o,k} &< i_{s,k} \\
j_{o,k} &> j_{s,k} \\
\endaligned\right\} \hskip.1in \text{with} \hskip.1in j_{o,k} - j_{s,k} < j_o
\hskip.1in \text{and} \hskip.1in j_{o,k} - b_k < j_o$$

$$(\gamma)_k \left\{\aligned
i_{o,k} &\leq i_{s,k} \\
j_{o,k} &\leq j_{s,k} \\
\endaligned \hskip.1in \& \hskip.1in 
\aligned
i_{o,k} &> a_k \\
j_{o,k} &< b_k \\
\endaligned
\right\} \hskip.1in \text{with} \hskip.1in i_{o,k} - a_k < i_o,$$

$$(\delta)_k \left\{\aligned
i_{o,k} &\leq i_{s,k} \\
j_{o,k} &\leq j_{s,k} \\
\endaligned \hskip.1in \& \hskip.1in 
\aligned
i_{o,k} &< a_k \\
j_{o,k} &> b_k \\
\endaligned
\right\} \hskip.1in \text{with} \hskip.1in j_{o,k} - b_k < j_o.$$

Then, $E_{p_k}$ being the exceptional divisor for $\roman{cp}_{1,k+1}$, obtained by
blowing up $p_k$, we have
$$\nu_{E_{p_k}}(R_{\log,f_1}) < \max\{i_o,j_o\}.$$
Moreover, let $p_{k+1} \in D_{X_{1,k+1}} \subset X_{1,k+1}$ be a point over $p_k$
where $g_{k+1}$ is in Subcase $2_{p_{k+1}}1_q2$.  Assume $g_{k+1}^{-1}(y_1,y_2)
\cdot \widehat{{\Cal O}_{X_{1,k+1},p_{k+1}}}$ is not principal.  Then we have
the conditions
$(i)_{k+1}$ with coordinate expression $(*)_{k+1}$ and
$(ii)_{k+1}$ (derived inductively from the conditions $(i)_k$ with coordinate
expression $(*)_k$ and
$(ii)_k$).

(Note that in the special case where in $g_0^*y_2$ there is no term
$x_{1,0}^ix_{2,0}^j$ with \linebreak
$\det\left[\matrix
a_0 & b_0 \\
i & j \\
\endmatrix\right] = 0$ and hence where subsequently in $g_k^*y_2$ there is no
term
$x_{1,k}^ix_{2,k}^j$ with
$\det\left[\matrix
a_k & b_k \\
i & j \\
\endmatrix\right] = 0$, we set $i_{s,0} = i_{s,k} = j_{s,0} = j_{s,k} =
\infty$ by convention.)
\endproclaim 

\proclaim{Diagram 3.3.6.C.3}\endproclaim

\vskip3in

\newpage

\demo{Proof}\enddemo We compute the coefficient
$\nu_{E_{p_k}}(R_{\log,f_1})$ of (the strict transform of) the exceptional
divisor $E_{p_k}$ for $\roman{cp}_{1,k+1}$, obtained by blowing up $p_k$, in
the logarithmic ramification divisor $R_{\log,f_1}$.  (We note that the
computation is identical to the one given for
$\nu_{E_p}(R_{\log,f_1}) =
\nu_{E_{p_0}}(R_{\log,f_1})$, except that at the end we compare
$\nu_{E_{p_k}}(R_{\log,f_1})$ to $\max\{i_o,j_o\}$ but not to
$\max\{i_{o,k},j_{o,k}\}$.).

\vskip.1in

Subcase $a_k + b_k \geq \min\{i_{o,k} + j_{o,k}, i_{s,k} + j_{s,k}\}$: In
this subcase, we have
$$\left\{\aligned
\nu_{E_{p_k}}(\frac{y_1}{y_2}) &= a_k + b_k - \min\{i_{o,k} + j_{o,k}, i_{s,k} +
j_{s,k}\}
\geq 0
\\
\nu_{E_{p_k}}(y_2) &= \min\{i_{o,k} + j_{o,k}, i_{s,k} + j_{s,k}\} > 0, \\
\endaligned\right.$$
which implies that the generic point of (the strict transform of) $E_{p_k}$ maps
under $f_1$ into the standard neighborhood of $q_1$ with a system of
regular parameters $(\frac{y_1}{y_2},y_2)$.  Therefore, we compute
$$\align
\nu_{E_{p_k}}(R_{\log,f_1}) &=
\nu_{E_{p_k}}\left(\frac{d(\frac{y_1}{y_2})}{(\frac{y_1}{y_2})} \wedge
\frac{dy_2}{y_2}/\frac{d(\frac{x_{1,k}}{x_{2,k}})}{(\frac{x_{1,k}}{x_{2,k}})} \wedge
\frac{dx_{2,k}}{x_{2,k}}\right) \\
&= \nu_{E_{p_k}}\left(\frac{dy_1}{y_1} \wedge
dy_2/\frac{dx_{1,k}}{x_{1,k}} \wedge \frac{dx_{2,k}}{x_{2,k}}\right) -
\nu_{E_{p_k}}(y_2)
\\ &= i_{o,k} + j_{o,k} - \min\{i_{o,k} + j_{o,k},i_{s,k} + j_{s,k}\}.\\
\endalign$$
Now we claim that
$$i_{o,k} + j_{o,k} - \min\{i_{o,k} + j_{o,k},i_{s,k} + j_{s,k}\} < \max\{i_o,j_o\}$$
unless $g_k^{-1}(y_1,y_2) \cdot \widehat{{\Cal O}_{X_{1,k},p_k}}$ is already
principal.

In fact, if
$$i_{o,k} + j_{o,k} \leq i_{s,k} + j_{s,k},$$
then
$$i_{o,k} + j_{o,k} - \min\{i_{o,k} + j_{o,k}, i_{s,k} + j_{s,k}\} = 0 <
\max\{i_o,j_o\}.$$ Therefore, we may assume
$$i_{o,k} + j_{o,k} > i_{s,k} + j_{s,k},$$
in which case we have
$$\align
i_{o,k} + j_{o,k} &- \min\{i_{o,k} + j_{o,k}, i_{s,k} + j_{s,k}\} = (i_{o,k} + j_{o,k})
- (i_{s,k} + j_{s,k})
\\ &= (i_{o,k} - i_{s,k}) + (j_{o,k} - j_{s,k}) \\
&=\left\{\aligned
&(i_{o,k} - i_{s,k}) - (j_{s,k} - j_{o,k}) < i_o \text{\ if\ }j_{s,k} \geq j_{o,k} \\
&(j_{o,k} - j_{s,k}) - (i_{s,k} - i_{o,k}) < j_o \text{\ if\ }i_{s,k} \geq i_{o,k} \\
& \text{or} \\
&g_k^{-1}(y_1,y_2) \cdot \widehat{{\Cal O}_{X_{1,k},p_k}} =
(x_{1,k}^{a_k}x_{2,k}^{b_k},x_{1,k}^{i_{s,k}}x_{2,k}^{j_{s,k}}) =
(x_{1,k}^{\min\{a_k,i_{s,k}\}}x_{2,k}^{\min\{b_k,j_{s,k}\}}) \\
&\hskip1.27in \text{principal\ }
\text{\ if\ }i_{s,k} < i_{o,k}
\hskip.05in \& \hskip.05in j_{s,k} < j_{o,k}.\\
\endaligned\right.\\
\endalign$$

Remark that, under the assumption
$$i_{o,k} + j_{o,k} > i_{s,k} + j_{s,k},$$
when $g_k^{-1}(y_1,y_2) \cdot \widehat{{\Cal O}_{X_{1,k},p_k}}$ is not
principal, we have
$$\align
j_{s,k} \geq j_{o,k} &\Longrightarrow (\alpha)_k \hskip.1in (\text{and\ hence\ }
i_{o,k} - i_{s,k} < i_o) \\ 
i_{s,k} \geq i_{o,k} &\Longrightarrow
(\beta)_k \hskip.1in (\text{and\ hence\ } j_{o,k} - j_{s,k} < j_o). \\
\endalign$$

\vskip.1in

Subcase $a_k + b_k < \min\{i_{o,k} + j_{o,k}, i_{s,k} + j_{s,k}\}$: In this
subcase, we have
$$\left\{\aligned
\nu_{E_{p_k}}(y_1) &= a_k + b_k > 0 \\
\nu_{E_{p_k}}(\frac{y_2}{y_1}) &= \min\{i_{o,k} + j_{o,k}, i_{s,k} + j_{s,k}\} - (a_k +
b_k) > 0,
\\
\endaligned\right.$$
which implies that the generic point of (the strict transform of) $E_{p_k}$ maps
under $f_1$ into the standard neighborhood of $q_1'$ with a system of
regular parameters $(y_1,\frac{y_2}{y_1})$.  Therefore, we compute
$$\align
\nu_{E_{p_k}}(R_{\log,f_1}) &=
\nu_{E_{p_k}}\left(\frac{dy_1}{y_1} \wedge
d(\frac{y_2}{y_1})/\frac{d(\frac{x_{1,k}}{x_{2,k}})}{(\frac{x_{1,k}}{x_{2,k}})} \wedge
\frac{dx_{2,k}}{x_{2,k}}\right) \\
&= \nu_{E_{p_k}}\left(\frac{dy_1}{y_1} \wedge
dy_2/\frac{dx_{1,k}}{x_{1,k}} \wedge \frac{dx_{2,k}}{x_{2,k}}\right) -
\nu_{E_{p_k}}(y_1)
\\ &= i_{o,k} + j_{o,k} - (a_k + b_k).\\
\endalign$$

Now we claim that
$$i_{o,k} + j_{o,k} - (a_k + b_k) < \max\{i_o,j_o\}$$
unless $g_k^{-1}(y_1,y_2) \cdot \widehat{{\Cal O}_{X_{1,k},p_k}}$ is already
principal.

In fact, we have
$$\align
i_{o,k} + j_{o,k} &- (a_k + b_k) = (i_{o,k} - a_k) + (j_{o,k} - b_k) \\
&=\left\{\aligned
&(i_{o,k} - a_k) - (b_k - j_{o,k}) < i_o \hskip.2in \text{\ if\ }b_k \geq j_{o,k} \\
&(j_{o,k} - b_k) - (a_k - i_{o,k}) < j_o \hskip.2in \text{\ if\ }a_k \geq i_{o,k} \\
&\text{or}\\
&g_k^{-1}(y_1,y_2) \cdot \widehat{{\Cal O}_{X_{1,k},p_k}} =
(x_{1,k}^{a_k}x_{2,k}^{b_k})
\\
& \hskip1.3in \text{\ principal\ }
\text{\ if\ }a_k < i_{o,k}
\hskip.05in \& \hskip.05in b_k < j_{o,k}\\
&(\text{Note\ that\ } a_k < i_{s,k} \hskip.05in \& \hskip.05in b_k < j_{s,k} \\
&\hskip.1in  \text{by\ the\ subcase\ assumption\ }a_k + b_k < \min\{i_{o,k} + j_{o,k},
i_{s,k} + j_{s,k}\}.)\\
\endaligned\right.\\
\endalign$$
Remark that, under the subcase assumption
$$a_k + b_k < \min\{i_{o,k} + j_{o,k}, i_{s,k} + j_{s,k}\},$$
when $g_k^{-1}(y_1,y_2) \cdot \widehat{{\Cal O}_{X_{1,k},p_k}}$ is not
principal, we have
$$\align
b_k \geq j_{o,k} &\Longrightarrow (\gamma)_k \hskip.1in (\text{and\ hence\ }i_{o,k} -
a_k < i_o) \\ 
a_k \geq i_{o,k} &\Longrightarrow
(\delta)_k \hskip.1in (\text{and\ hence\ } j_{o,k} - b_k < j_o). \\
\endalign$$
Therefore, we conclude
$$\nu_{E_{p_k}}(R_{\log,f_1}) < \max\{i_o,j_o\}.$$
Finally, it remains to prove the assertions $(i)_{k+1}$ with $(*)_{k+1}$ and
$(ii)_{k+1}$ on
\linebreak
$p_{k+1}
\in D_{X_{1,k+1}}
\subset X_{1,k+1}$, which is a point over $p_k$ and where $g_{k+1}$ is in Subcase
$2_{p_{k+1}}1_q2$.  When we blow up $p_k$, there appear two points (lying over $p_k$)
where $g_{k+1}$ is in Subcase $2_*1_q2$, one having the standard neighborhood with a
system of regular coordinates $(\frac{x_{1,k}}{x_{2,k}},x_{2,k})$ and the other having
the standard neighborhood with a system of regular coordinates
$(x_{1,k},\frac{x_{2,k}}{x_{1,k}})$.  By symmetry, we have only to check the
assertions on $p_{k+1}$ having the standard neighborhood with a system of
regular coordinates $(x_{1,k+1},x_{2,k+1}) = (\frac{x_{1,k}}{x_{2,k}},x_{2,k})$.

Firstly, by substituting
$$\left\{\aligned
&x_{1,k} = x_{1,k+1}x_{2,k+1} \\
&x_{2,k} = x_{2,k+1} \\
\endaligned\right.$$
into $(*)_k$, it is clear that the condition $(i)_{k+1}$ is satisfied with the system
of regular coordinates $(x_{1,k+1},x_{2,k+1})$ providing the coordinate expression
$(*)_{k+1}$.

Secondly, we check the condition $(ii)_{k+1}$, assuming that $g_{k+1}^{-1}(y_1,y_2)
\cdot \widehat{{\Cal O}_{X_{1,k+1},p_{k+1}}}$ is not principal.

Observe that
$$\align
a_{k+1} &= a_k,\\
b_{k+1} &= a_k + b_k,\\
i_{o,k+1} &= i_{o,k},\\
j_{o,k+1} &= i_{o,k} + j_{o,k}, \\
i_{s,k+1} &= i_{s,k},\\
j_{s,k+1} &= i_{s,k} + j_{s,k}. \\
\endalign$$

We analyze what happens at the $(k+1)$-the stage in each of the possibilities
$(\alpha)_k, (\beta)_k, (\gamma)_k, (\delta)_k$ at the $k$-th stage.

\vskip.1in

\noindent $(\alpha)_k$: Suppose that we are in the case $(\alpha)_k$.  

\vskip.1in

If $j_{o,k+1} < j_{s,k+1}$, then, since $i_{o,k+1} = i_{o,k} > i_{s,k} = i_{s,k+1}$, we
are in the case $(\alpha)_{k+1}$.  We also have
$$\align
i_{o,k+1} - i_{s,k+1} &= i_{o,k} - i_{s,k} < i_o \\
i_{o,k+1} - a_{k+1} &= i_{o,k} - a_k < i_o.\\
\endalign$$
If $j_{o,k+1} \geq j_{s,k+1}$, then, since $i_{o,k+1} = i_{o,k} > i_{s,k} = i_{s,k+1}$,
we have 
$$g_{k+1}^*y_2 = v_{k+1} \cdot x_{1,k+1}^{i_{s,k+1}}x_{2,k+1}^{j_{s,k+1}} \hskip.1in
\text{with} \hskip.1in v_{k+1} \in \widehat{{\Cal
O}_{X_{1,k+1},p_{k+1}}}^{\times} \hskip.1in
\text{unit}.$$   Therefore, the ideal
$$\align
g_{k+1}^{-1}(y_1,y_2) \cdot \widehat{{\Cal O}_{X_{1,k+1},p_{k+1}}} &=
(x_{1,k+1}^{a_{k+1}}x_{2,k+1}^{b_{k+1}}, x_{1,k+1}^{i_{s,k+1}}x_{2,k+1}^{j_{s,k+1}})
\\
&= (x_{1,k+1}^{\min\{a_{k+1},i_{s,k+1}\}}x_{2,k+1}^{\min\{b_{k+1},j_{s,k+1}\}}) \\
\endalign$$
is principal.

\vskip.1in
 
\noindent $(\beta)_k$: Suppose that we are in the case $(\beta)_k$.  

\vskip.1in

If $j_{o,k+1} > j_{s,k+1}$, then since $i_{o,k+1} = i_{o,k} < i_{s,k} = i_{s,k+1}$ we
are in the case $(\beta)_{k+1}$.  We also have
$$\align
j_{o,k+1} - j_{s,k+1} &= (i_{o,k} + j_{o,k}) - (i_{s,k} + j_{s,k}) \\
&= (j_{o,k} - j_{s,k}) - (i_{s,k} - i_{o,k}) < (j_{o,k} - j_{s,k}) < j_o \\
j_{o,k+1} - b_{k+1} &= (i_{o,k} + j_{o,k}) - (a_k + b_k) \\
&= (j_{o,k} - b_k) - (a_k - i_{o,k}) < (j_{o,k} - b_k) < j_o.\\
\endalign$$
(Note that under the case $(\beta)_k$ the inequality $a_k \leq i_{o,k}$ would imply
that \linebreak
$g_k^{-1}(y_1,y_2) \cdot \widehat{{\Cal O}_{X_{1,k},p_k}} =
(x_{1,k}^{a_k}x_{2,k}^{b_k})$ is principal, being against the assumption. 
Thus we have $a_k > i_{o,k}$.)

If $j_{o,k+1} \leq j_{s,k+1}$, then since $i_{o,k+1} = i_{o,k} < i_{s,k} = i_{s,k+1}$
we are in the case $(\delta)_{k+1}$, if $g_{k+1}^{-1}(y_1,y_2) \cdot
\widehat{{\Cal O}_{X_{1,k+1},p_{k+1}}}$ is not principal.  We also have
$$\align
j_{o,k+1} - b_{k+1} &= (i_{o,k} + j_{o,k}) - (a_k + b_k) \\
&= (j_{o,k} - b_k) - (a_k - i_{o,k}) < (j_{o,k} - b_k) < j_o.\\
\endalign$$
(See the note above.)

\vskip.1in

$(\gamma)_k$: Suppose that we are in the case $(\gamma)_k$.

\vskip.1in

Since, under the case $(\gamma)_k$,
$$i_{o,k} \leq i_{s,k} \hskip.1in \& \hskip.1in j_{o,k} \leq j_{s,k},$$
we have
$$\align
i_{o,k+1} &= i_{o,k} \leq i_{s,k} = i_{s,k+1} \\
j_{o,k+1} &= i_{o,k} + j_{o,k} \leq i_{s,k} + j_{s,k} = j_{s,k+1}. \\
\endalign$$
Therefore, supposing $g_{k+1}^{-1}(y_1,y_2) \cdot \widehat{{\Cal
O}_{X_{1,k+1},p_{k+1}}}$ is not principal, we are in the case
$(\gamma)_{k+1}$.  We also have
$$i_{o,k+1} - a_{k+1} = i_{o,k} - a_k < i_o.$$

\vskip.1in

$(\delta)_k$: Suppose that we are in the case $(\delta)_k$.

\vskip.1in

Since, under the case $(\delta)_k$,
$$i_{o,k} \leq i_{s,k} \hskip.1in \& \hskip.1in j_{o,k} \leq j_{s,k},$$
we have
$$\align
i_{o,k+1} &= i_{o,k} \leq i_{s,k} = i_{s,k+1} \\
j_{o,k+1} &= i_{o,k} + j_{o,k} \leq i_{s,k} + j_{s,k} = j_{s,k+1}. \\
\endalign$$
Therefore, supposing $g_{k+1}^{-1}(y_1,y_2) \cdot \widehat{{\Cal
O}_{X_{1,k+1},p_{k+1}}}$ is not principal, we are in the case
$(\delta)_{k+1}$.  We also have
$$\align
j_{o,k+1} - b_{k+1} &= (i_{o,k} + j_{o,k}) - (a_k + b_k) \\
&= (j_{o,k} - b_k) - (a_k - i_{o,k}) < (j_{o,k} - b_k) < j_o.\\
\endalign$$

Therefore, Lemma 3.3.6.C.2 is proved.

\vskip.1in

Now observe that at $p_0 = p$ where $g_0 = f$ is in Subcase $2_{p_0}1_q2$,
conditions $(i)_0$ with coordinate expression $(*)_0$ and $(ii)_0$ clearly
hold.  Therefore, by Lemma 3.3.6.C.2, at any point $p_k$ (over $p_0 = p$)
where the morphism $g_k$ is in Subcase
$2_{p_k}1_q2$, and where the ideal $g_k^{-1}(y_1,y_2) \cdot \widehat{{\Cal
O}_{X_{1,k},p_k}}$ is not principal, we see that conditions $(i)_k$ with
coordinate expression $(*)_k$ and $(ii)_k$ clearly hold, and that
$$\nu_{E_{p_k}}(R_{\log,f_1}) < \max\{i_o,j_o\}.$$

This completes the proof of the second part of statement b).

\vskip.1in

This completes the local analysis of our algorithm 2.1 in the case where
$\dim X = \dim Y = 2$.

\newpage

\S 4. Termination of the algorithm in the case where $\dim X = \dim Y = 2$

\vskip.1in

4.1. $\bold{Purpose\ of\ this\ section.}$ \hskip.1in The purpose of this
section is, based upon the local analysis in \S 3, to observe
that the logarithmic ramification divisor \it decreases \rm in the process
of the algorithm, under some partial order which we introduce among the Weil
divisors (on possibly different ambient spaces) and with respect to which
the descending chain condition is satisfied.  Thus the algorithm must
terminate after finitely many steps, ending with the logarithmic ramification
divisor being zero, i.e., achieving toroidalization of the given morphism.

\vskip.1in

4.2. $\bold{Partial\ order\ among\ the\ Weil\ divisors.}$ \hskip.1in Let
$R_1 =
\Sigma a_iF_i$ and $R_2 = \Sigma b_jG_j$ be Weil divisors on varieties $X_1$
and $X_2$, respectively, where the $F_i$ and $G_j$ are distinct irreducible
components with the coefficients in decreasing order, i.e.,
$$\align
&a_1 \geq a_2 \geq \cdot\cdot\cdot \geq a_{s-1} \geq a_s \\
&b_1 \geq b_2 \geq \cdot\cdot\cdot \geq b_{t-1} \geq b_t. \\
\endalign$$
Then we define the partial order among the Weil divisors so that
$$R_1 \geq R_2 \hskip.1in \text{if} \hskip.1in (a_1, a_2, ... , a_{s-1}, a_s) \geq
(b_1, b_2, ... , b_{t-1}, b_t),$$ where the second inequality is given with respect to
the lexicographical order. 

It is clear that the set
$$\{(a_1, a_2, ... , a_{s-1}, a_s)\}$$
satisfies the descending chain condition (i.e., there is no infinite strictly
decreasing sequence) and that so does the set of Weil divisors (possibly on different
ambient spaces). 

\vskip.1in

4.3. $\bold{Termination\ of\ the\ algorithm.}$

\proclaim{Theorem 4.3.1} The algorithm 2.1 for
toroidalization of a morphism
$f:(U_X,X) \rightarrow (U_Y,Y)$ in the logarithmic category, where
$\dim X = \dim Y = 2$, terminates after finitely many steps.

More precisely, in the process of the algorithm
$$\CD
(U_X,X) = (U_{X_0},X_0) & \hskip.1in \overset{\roman{cp}_1}\to{\leftarrow}
\hskip.1in & (U_{X_0},X_0) &
\hskip.1in \overset{\roman{cp}_2}\to{\leftarrow} &
\hskip.1in \cdot\cdot\cdot \hskip.1in &
\overset{\roman{cp}_i}\to{\leftarrow} \hskip.1in  & (U_{X_i},X_i) &
\hskip.1in \leftarrow &
\hskip.1in \cdot\cdot\cdot
\\ @VV{f = f_0}V @VV{f_1}V @. @VV{f_i}V \\
(U_X,X) = (U_{X_0},X_0) & \hskip.1in \overset{\roman{bp}_1}\to{\leftarrow}
\hskip.1in & (U_{X_0},X_0) &
\hskip.2in \overset{\roman{bp}_2}\to{\leftarrow} &
\hskip.1in \cdot\cdot\cdot \hskip.1in &
\overset{\roman{bp}_i}\to{\leftarrow} \hskip.1in & (U_{X_i},X_i) &
\hskip.1in \leftarrow &
\hskip.1in \cdot\cdot\cdot,
\\
\endCD$$
(o) the logarithmic ramification divisor never increases (with respect to the order
introduced in 4.2), i.e.,
$$R_{\log,f_i} \geq R_{\log,f_{i+1}},$$
(i) we have no infinite sequence of consecutive blowups
$$\roman{bp}_j \circ \roman{bp}_{j+1} \circ \cdot\cdot\cdot \circ \roman{bp}_{j + i}
\circ
\cdot\cdot\cdot$$
where the centers are only points $q$ of type $2_q$,

(ii) every time we have a point $q$ of type $1_q$ in the center of $\roman{bp}_i$,
the logarithmic ramification divisor strctly decreases (with respect to the order
introduced in 4.2), i.e.,
$$R_{\log,f_i} > R_{\log,f_{i+1}}.$$
Therefore, by the descending chain condition on the partial order, the algorithm must
terminate after finitely many steps.
\endproclaim

\demo{Proof}\enddemo (o) In each of the subcases from 3.3.1 through 3.3.6, we
conclude in the subsection C that
$$\text{either\ }R_{\log,f_i} = R_{\log,f_{i+1}} \text{\ or\ } R_{\log,f_i} >
R_{\log,f_{i+1}}.$$
(i) We may assume $j = 0$ for the proof.  Suppose that at the $0$-th stage, the
center of
$\roman{bp}_1$ consists only of points $q$ of type $2_q$.  This means that
$f(R_{\log,f_0})$ consists only of points of type
$2_q$.  

Suppose $p \not\in R_{\log,f_0}$, i.e., $f_0$ is toroidal in a neighborhood of
$p$.  Then by 3.3.1.B, over a neighborhood of $p$, the morphism $f_1$ remains
toroidal.  That is to say, over this neighborhood of $p$, there is no intersection
with $R_{\log,f_1}$.

Suppose $p \in R_{\log,f_0}$.  Then by the local analysis in 3.3.1, we have
coordinate expression
$$\left\{\aligned
f^*y_1 &= u \cdot x_1^{a_p}x_2^{b_p} \\
f^*y_2 &= v \cdot x_1^{c_p}x_2^{d_p} \\
\endaligned \right. \text{\ where\ }\det\left[\matrix
a_p & b_p \\
c_p & d_p \\
\endmatrix\right] = 0 \text{\ with\ units\ }u,v \in \widehat{{\Cal
O}_{X_0,p}}^{\times},$$ 
for some systems of regular parameters $(x_1,x_2)$ and
$(y_1,y_2)$ of
$\widehat{{\Cal O}_{X_0,p}}$ and
$\widehat{{\Cal O}_{Y_0,p}}$, compatible with the logarithmic structures.  

Remark that the matrix
$$\left[\matrix
a_p & b_p \\
c_p & d_p \\
\endmatrix\right]$$
is independent of the choice of the systems of regular parameters, up to the column
change and row change, and that the set
$$\{\left[\matrix
a_p & b_p \\
c_p & d_p \\
\endmatrix\right];p \in R_{\log,f_0}\}$$
is a finite set.

By the local analysis 3.3.1.B, the canonical principalization $\roman{cp}_1$ is an
isomorphism in a neighborhood of
$p \in R_{\log,f_0}$.

\vskip.1in

We have two possibilities: either

$\circ$ $p \not\in R_{\log,f_1}$

\noindent or 

$\circ$ $p \in R_{\log,f_1}$ where $f_1$ is in Subcases $* \hskip.03in 2_{q_1}
\hskip.03in *$ at $p$ with the corresponding matrix $\left[\matrix a_p - c_p & b_p -
d_p
\\ c_p & d_p \\
\endmatrix\right]$ (under the convention $a_p \geq c_p$ and $b_p \geq d_p$), since we
assume that $f_1(R_{\log,f_1})$ consists of points $q_1$ of type $2_{q_1}$.

But from this it follows that the length of the sequence could at most be
$$\max\{a_p + b_p + c_p + d_p; p \in R_{\log,f_0}\},$$
proving statement (i).

\vskip.1in

(ii) Statement (ii) follows from the local analysis in the subsection C of each
of the subcases from 3.3.2. through 3.3.6.

\vskip.1in

Now termination of the algorithm after finitely many steps is an easy consequence
of statements (i) and (ii), combined with the descending chain condition on the set
of the Weil divisors under the order introduced in 4.2.

\vskip.1in

4.4. $\bold{Corollaries}$

\proclaim{Corollary 4.4.1 (Resolution of singularities of morphisms in the
logarithmic category in dimension 2)} Let
$f:(U_X,X)
\rightarrow (U_Y,Y)$ be a morphism in the logarithmic category of nonsingular toroidal
embeddings with $\dim X = 2$.  Then there exist sequences of blowups with permissible
centers $\pi_X:(U_{X'},X') \rightarrow (U_X,X)$ and
$\pi_Y:(U_{Y'},Y')
\rightarrow (U_Y,Y)$ such that the induced map $f':(U_{X'},X')  \rightarrow (U_Y',Y')$
is log smooth.
\endproclaim

\demo{Proof}\enddemo When $\dim Y = 0$, there is nothing to prove.  When $\dim Y =
1$, it follows from 2.2.1 that $f$ is already log smooth.  When $\dim Y = 2$, the
assertion is a consequence of Theorem 4.3.1.

\proclaim{Corollary 4.4.2 (Toroidalization in dimension 2)} Toroidalization
conjecture holds for a dominant morphism $f:X \rightarrow Y$ between nonsingular
varieties with
$\dim X = 2$.
\endproclaim

\demo{Proof}\enddemo This follows from 1.5 and Corollary 4.4.1.

\vskip.1in

4.5. $\bold{Further\ remarks}$. \hskip.1in

\vskip.1in

\proclaim{Corollary 4.5.1 (Strong factorization of birational maps in dimension 2)}
Strong factorization of birational maps in dimension 2 holds: Let $\varphi:X
\dashrightarrow Y$ be a proper birational map between nonsingular varieties in
dimension 2.  Then there exist a sequence of blowups with points as centers
$$X = X_0 \overset{\roman{\phi}_1}\to{\leftarrow}
X_1 \overset{\roman{\phi}_2}\to{\leftarrow}
\cdot\cdot\cdot \overset{\roman{\phi}_{l-1}}\to{\leftarrow}
X_{l-1} \overset{\roman{\phi}_l}\to{\leftarrow} X_l$$ 
and a sequene of
blowdowns
$$X_l = Y_m \overset{\roman{\psi}_m}\to{\rightarrow} Y_{m-1}
\overset{\roman{\psi}_{m-1}}\to{\rightarrow} \cdot\cdot\cdot
\overset{\roman{\psi}_2}\to{\rightarrow} Y_1
\overset{\roman{\psi}_1}\to{\rightarrow} Y_0 = Y$$
such that
$$\varphi = \psi_1 \circ \cdot\cdot\cdot \circ \psi_m \circ \phi_l^{-1} \circ
\cdot\cdot\cdot \circ \phi_1^{-1}.$$
\endproclaim

\demo{Proof}\enddemo The whole point of
mentioning this well-known fact is to emphasize the way we prove it, according to the
general strategy Step I, which we have now completed proving toroidalization in
dimension 2, and step II, which is a consequence of an easy combonatorial fact of
geometry of convex bodies in dimension 2 (See 0.4 and Abramovich-Matsuki-Rashid[5]). 
We remark that the existing proofs for toroidalization use the factorization theorem
of proper birational morphisms via Castelnuovo's contractibility of a
$(-1)$-curve.

\vskip.1in

4.4.2. \it With the assumption of properness on $f$.\rm \hskip.1in  

\vskip.1in

Suppose that in the setting of 2.1.1 for a morphism $f:(U_X,X) \rightarrow (U_Y,Y)$
in the logarithmic category, we put the properness assumption on the morphism $f$.  

\vskip.1in

Then we have much less possibilities in 3.3.5) Subcase
$2_p1_q1$ and 3.3.6) Subcase $2_p1_q2$: Let $U_q$ be an open neighborhood of $q$ with
a system of regular parameters $(y_1,y_2)$ with $H_1 = \{y_1 = 0\}$.  Then since $f$
is smooth over $U_q - D_Y = U_q - H_1$, Abhyankar's lemma tells us that the
normalization
$\mu:\widetilde{U_q} \rightarrow U_q$ of $U_q$ in the function field $k(X)$ has a
system of regular parameters $(\widetilde{y_1},\widetilde{y_2})$ such that $\mu$ can
be written in these coordinate systems
$$\left\{\aligned
\mu^*y_1 &= \widetilde{y_1}^e \\
\mu^*y_2 &= \widetilde{y_2}. \\
\endaligned\right. \hskip.1in \text{with} \hskip.1in b > 0.$$
Now the factorization theorem of proper birational morphisms in dimension 2 tells us
that $g:f^{-1}(U_q) \rightarrow \widetilde{U_q}$ is a sequence of blowups of points.

3.3.5) Subcase $2_p1_q1$: Since $G_2$ is not exceptional for $f$ and hence not for
$g$, we see that there exists a system of regular parameters $(x_1,x_2)$ of
$\widehat{{\Cal O}_{X,p}}$ such that
$$\left\{\aligned
g^*\widetilde{y_1} &= x_1^{a'}x_2 \\
g^*\widetilde{y_2} &= x_1 \\
\endaligned\right.$$
Combining the above two, we conclude
$$\left\{\aligned
f^*y_1 &= x_1^ax_2^b \\
f^*y_2 &= x_1 \\
\endaligned\right. \hskip.1in \text{with} \hskip.1in a = a'e, b = e.$$

3.3.6) Subcase $2_p1_q2$: Since both $G_1$ and $G_2$ are exceptional for $f$ and
hence for
$g$, we see that there exists a system of regular parameters $(x_1,x_2)$ of
$\widehat{{\Cal O}_{X,p}}$ such that we have either
$$\left\{\aligned
f^*y_1 &= x_1^ax_2^b \\
f^*y_2 &= x_1^cx_2^d \\
\endaligned\right. \hskip.1in \text{with} \hskip.1in \det\left[\matrix
a & b \\
c & d \\
\endmatrix\right] \neq 0$$
where
$$\align
& a \geq c \hskip.1in \& \hskip.1in b \geq d, \text{\ or}\\
& a \leq c \hskip.1in \& \hskip.1in b \leq d.\\
\endalign$$
or
$$\left\{\aligned
f^*y_1 &= u \cdot x_1^ax_2^b \\
f^*y_2 &= v \cdot x_1^cx_2^d \\
\endaligned\right. \hskip.1in \text{with} \hskip.1in \det\left[\matrix
a & b \\
c & d \\
\endmatrix\right] = 0 \hskip.1in \text{and} \hskip.1in u,v \in \widehat{{\Cal
O}_{X,p}}^{\times} \hskip.1in
\text{units}.$$ 
The conclusion is that, under the properness assumption on $f$, in
both 3.3.5) Subcase
$2_p1_q1$ and 3.3.6) Subcase $2_p1_q2$ the ideal $f^{-1}(y_1,y_2) \cdot
\widehat{{\Cal O}_{X,p}}$ is already principal, and hence that the canonical
principalization is an isomorphism in a neighborhood of $p$.

Under the properness assumption on $f$, no new irreducible components appear in the
logarithmic ramification divisor in the process.  Moreover, identifying the strict
transforms of the irreducible components in the logarithmic ramification divisors, we
conclude (via the \it usual \rm order among the Weil divisors) 

\vskip.1in

3.3.1) Subcase $* \hskip.03in 2_q \hskip.03in *$ $\rightarrow $ $R_{\log,f_0} =
R_{\log,f_1}$

3.3.2) Subcase $1_p1_q0$ $\rightarrow $ $R_{\log,f_0} =
R_{\log,f_1}$

3.3.3) Subcase $1_p1_q1$ $\rightarrow $ $R_{\log,f_0} >
R_{\log,f_1}$

3.3.4) Subcase $2_p1_q0$: Impossible

3.3.5) Subcase $2_p1_q1$ $\rightarrow$ $R_{\log,f_0} >
R_{\log,f_1}$

3.3.6) Subcase $2_p1_q2$ $\rightarrow$ $R_{\log,f_0} >
R_{\log,f_1}$.

\vskip.1in

Therefore, the proof of Theorem 4.3.1 (ii) becomes much simpler under the
properness assumption on $f$.

\newpage

$$\bold{REFERENCES}$$

\vskip.1in

[1] S. Abhyankar, \it Simultaneous resolutions for algebraic surfaces, \rm Amer. J. Math
$\bold{78}$ (1956), 761-790

[2] S. Abhyankar, \it On the valuations centered in a local domain, \rm Amer. J. Math
$\bold{78}$ (1956), 321-348

[3] D. Abramovich and A.J. de Jong, \it Smoothness, semistability, and toroidal
geometry, \rm J. Alg. Geom. $\bold{6}$ (1997), 789-801

[4] D. Abramovich and K. Karu, \it Weak semistable reduction in characteristic zero, \rm
Invent. Math. $\bold{139}$ (2000), no. 2, 241-273

[5] D. Abramovich,K. Matsuki and S. Rashid, \it A note on the factorization theorem of
toric birational maps after Morelli and its toroidal extension, \rm Tohoku Math. J. (2)
$\bold{51}$ (1999), no. 4, 489-537

[6] D. Abramovich, K. Karu, K. Matsuki and J. W{\l}odarczyk, \it Torification and
factorization of birational maps, \rm J. Amer. Math. Soc. $\bold{15}$ (2002), 531-572

[7] S. Akbulut and H. King, \it Topology of algebraic sets, \rm MSRI publications 25

[8] E. Bierstone and P. Milman, \it Canonical desingularization in characteristic zero
by blowing up the maximum strata of a local invariant, \rm Invent. Math. $\bold{128}$
(1997), 207-302

[9] C. Christensen, \it Strong domination/weak factorization of three dimensional
regular local rings, \rm Journal of the Indian Math. Soc. $\bold{45}$ (1981), 21-47

[10] S.D. Cutkosky, \it Local factorization of birational maps, \rm Advances in Math.
$\bold{132}$ (1997), 167-315

[11] S.D. Cutkosky, \it Local monomialization and factorization of morphisms, \rm
Ast\'erisque 260, Soc. Math. France (1999)

[12] S.D. Cutkosky, \it Monomialization of morphisms from 3-folds to surfaces, \rm
preprint math.AG/0010002

[13] S.D. Cutkosky and O. Piltant, \it Monomial resolutions of morphisms of algebraic
surfaces, \rm Special issue in honor of Robin Hartshorne, Comm. Algebra $\bold{28}$
(2000), no. 12, 5935-5959

[14] S. Encinas and O. Villamayor, \it A course on constructive desingularization and
equivariance, \rm in Resolution of singularities (Obergurgl, 1997), Progress in
Math. $\bold{181}$ (2000), Birkh\"auser, 147-227

[15] S. Encinas and O. Villamayor, \it A new theorem of desingularization over fields of
characteristic zero, \rm preprint (1999)

[16] H. Hironaka, \it Resolution of singularities of an algebraic variety over a field
of characteristic zero, \rm Annals of Math. $\bold{79}$ (1964), 109-324  

[17] S. Iitaka, \it Algebraic Geometry (An introduction to birational geometry of
algebraic varieties), \rm Graduate Texts in Math. 76  (1982), Springer

[18] A.J. de Jong, \it Smoothness, semistability, and alterations, \rm Publ. Math.
I.H.E.S. $\bold{83}$ (1996), 51-93

[19] K. Karu, \it Local strong factorization of birational maps, \rm preprint (2003)

[20] K. Kato, \it Toric singularities, \rm Amer. J. Math. $\bold{116}$ (1994), 1073-1099

[21] G. Kempf, F. Knudsen, D. Mumford and B. Saint-Donat, \it Toroidal embeddings I, \rm
Lecture Notes in Math. $\bold{339}$ (1973), Springer

[22] H. King, \it Resolving singularities of maps, \rm Real algebraic geometry and
topology (East Lansing, Michigan, 1993), Contemp. Math. (1995), Amer. Math. Soc.

[23] K. Matsuki, Correction to ``A note on the factorization theorem of birational maps
after Morelli and its toroidal extension", \rm Tohoku Math. J. $\bold{52}$ (2000),
629-631

[24] K. Matsuki, \it Lectures on factorization of birational maps, \rm RIMS preprint
$\bold{1281}$ (1999)

[25] K. Matsuki, \it Notes on the inductive algorithm of resolution of singularities by
S. Encinas and O. Villamayor, \rm preprint (2002)

[26] R. Morelli, \it The birational geometry of toric varieties, \rm J. Alg. Geom.
$\bold{5}$ (1996), 751-782

[27] T. Oda, \it Convex bodies and algebraic geometry, \rm Ergebnisse der Mathematik und
ihrer Grenzgebiete $\bold{15}$ (1988)

[28] J. W{\l}odarczyk, \it Decomposition of birational toric maps in blow-ups and
blow-downs (A proof of the weak Oda conjecture), \rm Transaction of the AMS
$\bold{349}$ (1997), 373-411

[29] J. Wlodarczyk, \it Birational cobordism and factorization of birational maps, \rm
J. Alg. Geom. $\bold{9}$ (2000), no. 3, 425-449

[30] J. W{\l}odarczyk, \it Toroidal varieties and the weak factorization theorem, \rm
preprint math.AG/9904076

[31] O. Zariski, \it Algebraic surfaces, \rm Springer-Verlag (1934)

\enddocument